\documentclass[smallextended,envcountsect]{styles/svjour3}
\usepackage[
    a4paper,
    left=2.5cm,
    right=2.5cm,
    top=2.5cm,
    bottom=2.5cm
]{geometry}
\smartqed
\usepackage[utf8]{inputenc}
\usepackage[english]{babel}

\usepackage{cite}
\usepackage{amsmath,amssymb,amsfonts}

\usepackage{amsthm}
\usepackage{xpatch}
\usepackage{dsfont}
\usepackage{acronym}
\usepackage{mathtools}
\usepackage{units}
\usepackage{graphicx}
\usepackage{algorithm}
\usepackage{algpseudocode}
\usepackage{hyperref}
\usepackage{cleveref}
\usepackage{url}            
\usepackage{booktabs}       
\usepackage{nicefrac}       
\usepackage{microtype}      
\usepackage{textcomp}
\usepackage[short,nocomma]{styles/optidef}
\usepackage{xcolor}
\definecolor{oxfordblue}{RGB}{0,33,71}          
\usepackage{tikz}
\usetikzlibrary{calc, plotmarks, positioning, fit, fpu, arrows.meta}
\usepackage{pgfplots}
\usepgfplotslibrary{statistics, fillbetween, groupplots, external, colorbrewer, colormaps}
\pgfplotsset{compat=1.18}
\usepackage{xparse}
\usepackage{booktabs}
\usepackage{enumitem}
\setlist[enumerate]{leftmargin=.5in}
\setlist[itemize]{leftmargin=.5in}

\DeclareRobustCommand\Measure[3]{\ifthenelse{\isempty{#2}}%
  {\ifthenelse{\isempty{#3}}
  {#1}
  {#1(#3)}}
  {\ifthenelse{\isempty{#3}}
  {#1_{#2}}
  {#1_{#2}({#3})}}
}
\renewcommand\P[2]{\Measure{\mathbb{P}}{#1}{#2}}
\newcommand\Q[2]{\Measure{\mathbb{Q}}{#1}{#2}}

\newcommand\PQ[2]{\Measure{\pi}{#1}{#2}}

\newcommand\dP[2]{\Measure{\,\mathrm{d}\mathbb{P}}{#1}{#2}}
\newcommand\dQ[2]{\Measure{\,\mathrm{d}\mathbb{Q}}{#1}{#2}}
\newcommand\dPQ[2]{\Measure{\,\mathrm{d}\pi}{#1}{#2}}

\newcommand\E[2]{{\mathbb{E}_{#1}\left[{#2}\right]}}
\newcommand\Eonly[1]{{\mathbb{E}_{#1}}}

\NewDocumentCommand{\Risk}{ m m o o}{
  \mathcal{R}_{#1}
  \IfNoValueTF{#3}{}{^{#3}}
  \big[{#2}
  \IfNoValueTF{#4}{\big]}{}
}
\newcommand\Unc[3]{
 \ifthenelse{\isempty{#3}}%
  {\mathcal{U}_{#1}\big[{#2}\big]}
  {\mathcal{U}_{#1}^{#3}\big[{#2}\big]}
  }

\newcommand\Opt[4]{
  \ifthenelse{\isempty{#2}}%
  {\mathrm{{#1}}}
  {\underset{\substack{\displaystyle #2}}{\mathrm{{#1}}}}
  \;
  \ifthenelse{\isempty{#4}}%
  {{#3}}
  {\Big\{ {#3}\  \Big| \  {#4} \Big\}}
}

\newcommand\transpose[1]{#1^\mathsf{T}}
\newcommand\Id[1]{\mathrm{I_{#1}}}
\newcommand\e[1]{\exp \left( #1 \right)}
\newcommand\KL[2]{\mathrm{KL} \left( #1 \left\| #2 \right. \right)}
\newcommand\MutualInf[3]{\mathrm{I}_{#1} \left( #2 ; #3 \right)}

\newcommand\ShannonEnt[2]{\mathrm{H}_{#1} \left( #2 \right)}
\newcommand\DiffEnt[2]{\mathrm{h}_{#1} \left( #2 \right)}
\renewcommand\d[1]{\,\mathrm{d}{#1}} 

\crefname{fact}{Fact}{Facts}
\crefname{equation}{Equation}{Equations}
\crefname{proposition}{Proposition}{Propositions}
\crefname{corollary}{Corollary}{Corollaries}
\crefname{remark}{Remark}{Remarks}
\crefname{hypothesis}{Hypothesis}{Hypotheses}
\crefname{assumption}{Assumption}{Assumptions}
\crefname{section}{Section}{Sections}
\crefname{subsection}{Section}{Sections}
\crefname{algorithm}{Algorithm}{Algorithms}
\crefname{theorem}{Theorem}{Theorems}
\crefname{figure}{Figure}{Figures}

\renewcommand{\makeheadbox}{}
\journalname{JOTA}

\newcommand\AMEPOptConst{\frac{d}{p}\,\operatorname{vol}\!\left(B_p^d\right)\,\Gamma\!\left(\frac{d}{p}\right)\,\left(\frac{p}{d}\right)^{\!d/p}\,\zeta^{d}}
\newcommand\AMEPOptConstInv{\frac{p}{d\,\operatorname{vol}(B_p^d)\,\Gamma(d/p)} \left(\frac{d}{p}\right)^{d/p}\zeta^{-d}}
\newcommand\AMEPOptVal{\frac{d}{p} + \log \left(\AMEPOptConst\right)}
\newcommand\AMEPOptArg{\AMEPOptConstInv\,\int_C \e{-\frac{d}{p\zeta^p} \|x-v\|^p_p} \d{v}}

\begin{document}
\title{Transport, Entropy, and the Cost of Dependence}
\subtitle{Maximum Conditional Entropy under Transport Constraints}

\author{Andr\'e Bertolace, Konstantinos Gatsis, Kostas Margellos}

\institute{Andr\'e Bertolace (andre.bertolace@eng.ox.ac.uk)
            \at
           Department of Engineering, University of Oxford,Oxford, OX1 3PJ, U.K.
           \and
           Konstantinos Gatsis \at
             Electrical and Computer Engineering, Villanova University, Villanova, PA, USA.
           \and
           Kostas Margellos \at
           Department of Informatics, National and Kapodistrian University of Athens, Panepistimioupolis Ilisia, Athens, 16122, and Department of Engineering, University of Oxford, Oxford, OX1 3PJ, U.K.
}
\date{October 2026}

\maketitle
\begin{abstract}
Maximum-entropy inference under Wasserstein constraints can be computationally challenging due to the underlying optimization over transport couplings. We investigate a tractable alternative based on maximizing conditional entropy under an average transport constraint. On Euclidean spaces, this formulation admits an explicit generalized-Gaussian additive-noise solution, providing a feasible approximation to the Wasserstein-constrained problem without solving an optimal transport problem. We establish quantitative bounds on the resulting optimality gap by showing that the two formulations admit a common feasible-marginal representation, with their discrepancy characterized by the mutual information retained between the nominal and perturbed variables. For discretely supported nominal distributions, we further obtain an exact decomposition of the gap into a marginal-entropy gain and a weighted Jensen--Shannon divergence.
\end{abstract}

\section{Introduction}
\label{sec:introduction}

Selecting a probability distribution from incomplete information is a fundamental problem in statistical inference and statistical physics. When the available information does not uniquely determine a distribution, the maximum entropy principle (MEP) provides a systematic criterion for selecting the least informative distribution consistent with that information \cite{jaynes1957part1,jaynes1957part2,Csiszar1975,shore1980,skilling1988,Book_CoverThomas2006}. Traditionally, the available information is expressed through moment or expectation constraints, leading to Gibbs-type distributions. In many settings, however, uncertainty is naturally described relative to a prescribed nominal distribution, with admissible distributions restricted to a neighborhood of that nominal law. Optimal transport provides a geometric framework for describing such neighborhoods through the Wasserstein distance, which quantifies the minimum cost required to transport one probability distribution into another \cite{Villani2008,Santambrogio2015}.

This motivates the study of maximum-entropy inference under Wasserstein constraints. Given a nominal distribution $\P{X}{}$, the objective is to identify a maximally uncertain distribution $\P{V}{}$ whose Wasserstein distance from $\P{X}{}$ does not exceed a prescribed radius. Such a formulation combines the maximum entropy principle with a geometric description of distributional uncertainty. However, directly solving the resulting optimization problem can be challenging. The Wasserstein constraint is defined through an optimization over transport couplings, and characterizing the entropy-maximizing distribution generally requires solving a variational problem involving the induced marginal distribution and the associated transport geometry. This raises a natural question: \emph{can Wasserstein-constrained maximum-entropy inference be approximated through a computationally simpler formulation, while retaining quantitative guarantees on the loss of optimality?}

One alternative is to impose the transport budget directly on the stochastic transition generating the perturbed variable. Given a nominal random variable $X$, one seeks a conditional distribution $\P{V|X}{}$ maximizing the conditional entropy of $V$ given $X$, subject to a bound on the expected transport cost between $X$ and $V$. This formulation is closely related to the maximum-caliber principle, originally proposed by Jaynes \cite{Jaynes1980}, which extends maximum-entropy reasoning to stochastic trajectories and dynamical laws subject to constraints on transition or path observables \cite{Stock2008,Ge2012,Presse2013}. We refer to this formulation as the average-transport maximum entropy principle (A-MEP). Unlike the Wasserstein constraint, which involves minimizing the transport cost over all couplings between two distributions, the average-transport constraint is evaluated directly under the joint law induced by the perturbation kernel. This distinction leads to a particularly tractable optimization problem, on Euclidean spaces, A-MEP admits an explicit generalized-Gaussian additive-noise solution, including Gaussian and Laplace perturbations as special cases.

In contrast, Wasserstein-constrained maximum-entropy inference imposes the transport budget only on the marginal distribution induced by the perturbation. This distribution-level perspective is related to Wasserstein variational principles, including the Jordan--Kinderlehrer--Otto scheme \cite{JKO1998}, in which probability distributions are selected through variational problems involving Wasserstein geometry. We refer to the corresponding formulation as the Wasserstein maximum entropy principle (W-MEP). Importantly, the coupling certifying Wasserstein proximity need not coincide with the joint law generated by the entropy-bearing transition kernel. Consequently, although A-MEP and W-MEP maximize the same conditional-entropy objective, they impose different restrictions on the stochastic relationship between the nominal and perturbed variables.

The computational simplicity of A-MEP makes it an attractive surrogate for W-MEP. Indeed, every perturbation kernel satisfying the average-transport constraint induces a marginal distribution lying within the corresponding Wasserstein neighborhood of the nominal law. Thus, the explicit A-MEP solution provides a feasible construction for W-MEP without requiring the computation of an optimal transport coupling or the solution of the Wasserstein-constrained variational problem. Nevertheless, the two formulations need not yield the same optimal value or the same entropy-maximizing marginal distribution. Feasibility alone is therefore insufficient to establish the quality of the approximation. The central question addressed in this paper is: \emph{how accurately can the tractable A-MEP formulation approximate W-MEP, and how can the resulting loss of optimality be quantified?}

We address this question through three main contributions:
\begin{enumerate}
    \item First, we characterize the optimal solutions of both formulations. On Euclidean spaces, A-MEP admits an explicit generalized-Gaussian additive-noise solution, providing a computationally tractable construction of a feasible distribution for W-MEP. For W-MEP, we derive a dual characterization in terms of a Kantorovich potential and show that the entropy-maximizing conditional kernel can be chosen independently of the nominal realization. 
    \item Second, we establish that both formulations generate the same set of attainable marginal distributions and admit representations over a common feasible set. This reveals that their discrepancy is governed by the mutual information retained between the nominal and perturbed variables. We use this characterization to derive bounds on the entropy gap, including an explicit upper bound that provides a quantitative optimality guarantee for the A-MEP surrogate.
    \item Third, for finitely supported nominal distributions, we obtain finite-mixture representations and an exact decomposition of the entropy gap into a marginal-entropy gain and a weighted Jensen--Shannon divergence. A symmetric two-point example illustrates that the gap can be strict and that the formulations can yield different entropy-maximizing distributions.
\end{enumerate}.  

Related connections between entropy and optimal transport arise in regularized optimal transport and the Schrödinger problem \cite{Cuturi2013,Peyre2019,Leonard2014}, while entropy-based selection over Wasserstein balls has also been studied through minimum-cross-entropy formulations \cite{Vargas2021}. These works differ from the question considered here. Rather than introducing an entropy regularization of the transport problem or studying relative-entropy projections onto Wasserstein neighborhoods, we investigate whether a maximum-entropy problem under an average transport constraint can serve as a tractable approximation to its Wasserstein-constrained counterpart. Our results establish that this approximation can be characterized and bounded through information-theoretic quantities, providing a connection between computational tractability, transport geometry, and the dependence retained by stochastic perturbations.
\section{Problem setting}
\label{section:problem_formulation}

Throughout the paper, all random objects are defined on a common probability space $(\Omega,\mathcal{F},\P{}{})$. Let $(\mathcal{X}, \| \cdot \|)$ be a Polish space with Borel $\sigma$-algebra $\mathcal{B}(\mathcal{X})$, and let $X,V:(\Omega,\mathcal{F})\to(\mathcal{X},\mathcal{B}(\mathcal{X}))$ be random variables. Denote by $\mathcal{P}(\mathcal{X})$ the set of all Borel probability measures on $(\mathcal{X},\mathcal{B}(\mathcal{X}))$, and by $\mathcal{P}_p(\mathcal{X})$ the set of probability measures on $\mathcal{X}$ with finite $p$-th moment,
\[
\mathcal{P}_p(\mathcal{X}) = \left\{ \P{}{} \in \mathcal{P}(\mathcal{X})
: \int_{\mathcal{X}} \|x\|^p \dP{}{x} < \infty \right\}.
\]
For a given a nominal distribution $\P{X}{} \in \mathcal{P}_p(\mathcal{X})$, a perturbation is described by a stochastic kernel $\P{V|X}{}$ from $\mathcal{X}$ to $\mathcal{X}$. Together, $\P{X}{}$ and $\P{V|X}{}$ induce the joint law
\[
    \P{X,V}{} := \P{X}{} \otimes \P{V|X}{},
\]
defined by
\[
    \P{X,V}{B \times C} = \int_B \P{V|X=x}{C} \dP{X}{x}, \qquad B,C \in \mathcal{B}(\mathcal{X}),
\]
and its marginal with respect to $V$, denoted by $\P{V}{}$, satisfies
\[
    \P{V}{C} = \int_{\mathcal{X}} \P{V|X=x}{C} \dP{X}{x}, \qquad C \in \mathcal{B}(\mathcal{X}).
\]
We write
\[
\mathcal{K}_p(\mathcal{X}) = \left\{ \P{V|X}{} \text{ stochastic kernels on } \mathcal{X} : \int_{\mathcal{X}} \int_{\mathcal{X}} \|v\|^p \dP{V|X}{v} \dP{X}{x} < \infty \right\},
\]
for the set of perturbation kernels whose induced marginal belongs to
$\mathcal{P}_p(\mathcal{X})$.

The $p$-Wasserstein distance between probability measures $\P{X}{},\P{V}{}\in\mathcal P_p(\mathcal X)$ is defined by
\begin{customopti*}|s|
    {inf}
    {\PQ{}{} \in \Pi(\P{X}{}, \P{V}{})}
    {\int_{\mathcal{X}\times\mathcal{X}}
    \|x-v\|_p^p\,\dPQ{}{x,v}}
    {}
    {W_p(\P{X}{}, \P{V}{}) = \Big(}
    {\Big)^{1/p}}
    {},
\end{customopti*}
where $\Pi(\P{X}{},\P{V}{})$ denotes the set of all couplings of $\P{X}{}$ and $\P{V}{}$.

A distinction used throughout the paper is that between the perturbation-induced joint law $\P{X,V}{} = \P{X}{} \otimes \P{V|X}{}$ and a generic transport coupling of the same marginals $\pi \in \Pi(\P{X}{},\P{V}{})$. Since $\mathcal{X}$ is a standard Borel space, such a coupling admits a disintegration \cite[Thm. 8.5]{Kallenberg2021} of the form
\[
    \dPQ{}{x,v} = \dP{X}{x} \dQ{V|X=x}{v},
\]
for some stochastic kernel $\Q{V|X}{}$. In general, this $\Q{V|X}{}$ need not coincide with the perturbation kernel $\P{V|X}{}$ even though both kernels induce the same marginal $\P{V}{}$. Thus, $\P{X,V}{}$ is the coupling generated by the perturbation, whereas $\pi$ is a generic coupling used in the transport problem.

\paragraph{Entropy and mutual information.}
Entropy is not defined on an arbitrary Polish space without specifying additional structure. We therefore follow the standard information-theoretic convention of using $\ShannonEnt{\P{X,V}{}}{V}$ for Shannon entropy on finite or countable spaces and $\DiffEnt{\P{X,V}{}}{V}$ for differential entropy on Euclidean spaces \cite{Book_CoverThomas2006}. Mutual information is denoted by $\MutualInf{\P{X,V}{}}{X}{V}$ in both settings. Throughout, the subscript identifies the joint probability law under which the corresponding information-theoretic quantity is evaluated; marginal and conditional laws are understood to be induced by this joint law. Thus, for example, $\DiffEnt{\P{X,V}{}}{V}$ denotes the differential entropy of the $V$-marginal of $\P{X,V}{}$, while $\DiffEnt{\P{X,V}{}}{V|X}$ denotes the corresponding conditional differential entropy. For instance, if $\mathcal{X}$ is finite or countable, the Shannon entropy of $V \sim \P{V}{}$ is
\[
    \ShannonEnt{\P{X,V}{}}{V} := -\sum_{v\in\mathcal{X}} \P{V}{v}\log \P{V}{v}.
\]
If instead $\mathcal{X}\subseteq\mathbb{R}^d$ and $\P{V}{}$ admits a density $p_V$, its differential entropy is
\[
    \DiffEnt{\P{X,V}{}}{V} := -\int_{\mathcal{X}} p_V(v) \log p_V(v) \d{v},
\]
whenever the integral is well defined. Analogously, we also define the conditional entropy, in the discrete setting, 
\[
\ShannonEnt{\P{X,V}{}}{V | X} = - \sum_{x,v} \P{X,V}{x,v} \log \P{V|X=x}{v},
\]
as well as in the Euclidean setting, if $\P{V|X=x}{}$ admits a density $p_{V|X=x}(v)$ for $\P{X}{}$-almost every $x$, 
\[
\DiffEnt{\P{X,V}{}}{V | X} = - \int_{\mathcal{X}\times\mathcal{X}} \log p_{V|X=x}(v) \dP{X,V}{x,v}.
\]
Furthermore, the mutual information is
\[
    \MutualInf{\P{X,V}{}}{X}{V} := \KL{\P{X,V}{}}{\P{X}{}\otimes \P{V}{}},
\]
with the Kullback--Leibler divergence, 
\[
\KL{\P{}{}}{\Q{}{}} =
\begin{cases}
\displaystyle \int_\Omega \log \left(\frac{\dP{}{}}{\dQ{}{}}\right), \dP{}{},
& \text{if } \P{}{} \ll \Q{}{} \text{ and } \log \left(\frac{\dP{}{}}{\dQ{}{}}\right) \in L^1(\P{}{}), \\
\infty, & \text{otherwise.}
\end{cases}
\]
Whenever the relevant quantities are well defined, the entropy decomposition
\[
    \DiffEnt{\P{X,V}{}}{V | X} = \DiffEnt{\P{X,V}{}}{V}-\MutualInf{\P{X,V}{}}{X}{V},
\]
holds in the discrete and continuous Euclidean setting, respectively. 

\subsection{Problem setup}
We now introduce the two maximum-entropy problems studied in this paper.

\paragraph{Maximum entropy under average transport constraint (A-MEP).}
Under the average-transport formulation, both entropy and transport cost are evaluated under the joint law generated by the perturbation kernel, 
\begin{customopti}|s|
    {sup}
    {\P{V|X}{} \in \mathcal{K}_p(\mathcal{X})}
    {\DiffEnt{\P{X,V}{}}{V|X}}
    {}
    {A_\zeta := }
    {\label{equation:average_constrained_entropy}}
    \addConstraint{\E{X,V}{\|x-v\|_p^p} \leq \zeta^p.}
\end{customopti}

\paragraph{Maximum entropy under Wasserstein constraint (W-MEP).}
Under the Wasserstein formulation, the problem maximizes the same conditional entropy, but constrains only the marginal distribution induced by the perturbation kernel,
\begin{customopti}|s|
    {sup}
    {\P{V|X}{} \in \mathcal{K}_p(\mathcal{X})}
    {\DiffEnt{\P{X,V}{}}{V|X}}
    {}
    {W_\zeta :=}
    {\label{equation:wasserstein_constrained_entropy}}
    \addConstraint{W_p(\P{X}{}, \P{V}{}) \leq \zeta.}
\end{customopti} 

Note that, every kernel feasible for A-MEP induces a marginal feasible for W-MEP, since the joint law $\P{X}{}\otimes\P{V|X}{}$ is an admissible coupling between $\P{X}{}$ and $\P{V}{}$. Consequently,
\[
    W_p^p(\P{X}{},\P{V}{}) \leq \E{X,V}{\|x-v\|_p^p} \leq \zeta^p.
\]
Thus, the feasible set of A-MEP is contained in that of W-MEP. The inclusion can be strict because the joint law induced by $\P{V|X}{}$ need not be an optimal coupling between its marginals.

\section{Maximum entropy principle under transport constraints}
\label{section:mep}

\subsection{Average transport constrained MEP}
\label{subsection:avergate_mep}

In the continuous setting, the A-MEP problem admits a closed-form solution. The dual variational formulation provides the structure needed to identify the optimizer, and the next theorem shows that the entropy-maximizing kernel on $\mathbb{R}^d$ is of Gibbs type.

\begin{proposition}[A-MEP]
\label{proposition:AMEP}
Consider $\mathcal{X} = \mathbb{R}^d$, then, given $\P{X}{} \in \mathcal{P}_p(\mathcal{X})$, the optimization problem in \eqref{equation:average_constrained_entropy} admits an optimal value equal to,
\begin{equation*}
A_\zeta = \AMEPOptVal,
\end{equation*}
 with optimal solution given by a generalized gaussian distribution,
\begin{equation*}
\P{V|X}{}^\emph{A}(C) = \AMEPOptArg,
\end{equation*}
where $\Gamma(d) = \int_{0}^{\infty} x^{d-1} e^{-x} \d{x}  $ is the gamma function and $
\mathrm{vol}(B_p^d)=\frac{\left(2\Gamma\left(1+\frac{1}{p}\right)\right)^d}{\Gamma\left(1+\frac{d}{p}\right)}$ is the volume of the $\ell_p$-unit-ball.
\end{proposition}
The proof is provided in \cref{proof:proposition:AMEP}. 

The next corollaries show special cases of \cref{proposition:AMEP} corresponding to $p=2$, $p=1$ and $p \rightarrow \infty$.

\begin{corollary}[$\ell_2$ A-MEP]
\label{corollary:AMEP2}
Consider $\mathcal{X} = \mathbb{R}^d$, then, given $\P{X}{}  \in \mathcal{P}_2(\mathcal{X})$, the optimization problem in \eqref{equation:average_constrained_entropy}, with $p=2$, admits an optimal value $\frac{d}{2} + d\,\log \Big(\sqrt{\frac{2\pi}{d}} \zeta\Big)$,
and optimal solution given by $\P{V|X}{}^\emph{A} \sim \mathcal N\!\left(x,\frac{\zeta^2}{d}\Id{d}\right)$.
\end{corollary}

\begin{corollary}[$\ell_1$ A-MEP]
\label{corollary:AMEP1}
Consider $\mathcal{X} = \mathbb{R}^d$, then, given $\P{X}{}  \in \mathcal{P}_1(\mathcal{X})$, the optimization problem in \eqref{equation:average_constrained_entropy}, with $p=1$, admits an optimal value $d + d\,\log \frac{2\zeta}{d}$, an optimal optimal solution given by
$\P{V|X}{}^\emph{A} \sim \mathrm{Laplace}_d\left(x,\frac{\zeta}{d}\right)$,
equivalently, $\P{v_i|x_i}{}^\emph{A} \stackrel{\mathrm{i.i.d.}}{\sim} \mathrm{Laplace}\!\left(x_i,\frac{\zeta}{d}\right)$.
\end{corollary}

\begin{corollary}[Hard $\ell_\infty$ MEP]
\label{corollary:AMEPinfty}
Let $\zeta>0$ and let $\P{X}{}$ be a Borel probability measure on $\mathbb{R}^d$. The hard-constrained maximum conditional entropy problem
\begin{customopti*}|s|
    {sup}
    {\P{V|X}{} \in \mathcal{K}_p(\mathcal{X})}
    {\DiffEnt{\P{X,V}{}}{V|X}}
    {}
    {}
    \addConstraint{\|X-V\|_\infty \leq \zeta \qquad \P{X,V}{}\text{-a.s.}},
\end{customopti*}
admits the unique optimal conditional density, $V|X=x \sim \operatorname{Unif}\bigl(x+[-\zeta,\zeta]^d\bigr)$, and the corresponding optimal value is $d\log(2\zeta)$.
\end{corollary}

These results show that the entropy-maximizing conditional distribution under an average transport constraint admits a simple additive-noise interpretation. For the squared average transport cost, the maximizer is a multivariate Gaussian distribution with diagonal covariance, assigning variance $\zeta^2/d$ to each coordinate. For the absolute average transport cost, the optimal conditional law is an independent Laplace distribution with scale parameter $\zeta/d$ in each coordinate. The proofs are in \cref{proof:corollary:AMEP2,proof:corollary:AMEP1,proof:corollary:AMEPinfty}, respectively. 

\subsection{Wasserstein constrained MEP}
\label{section:mep_wasserstein}

The marginal Wasserstein problem no longer admits a purely local Gibbs kernel depending only on the displacement $x-v$. This is natural because, unlike the average-transport constraint, the Wasserstein constraint is imposed at the level of the induced marginal rather than on each transition from $x$ to $v$. Consequently, the conditional laws associated with different realizations of $X$ are coupled through the common marginal they induce, so the entropy--transport trade-off can no longer be resolved pointwise in $x$. Instead, the optimizer is characterized through a Kantorovich-type dual potential. The next result makes this distinction explicit by deriving a dual formulation in which the entropy contribution appears through a $\log$-partition function, while the transport contribution appears through a $c$-transform.

In this sense, the average transport formulation identifies the entropy kernel with the transport kernel, whereas the marginal Wasserstein formulation only requires their induced marginals to agree. The latter is particularly natural when individual transitions cannot be observed or matched, or when the available information concerns only population-level distributions rather than individual perturbations.

\begin{proposition}[W-MEP dual characterization]
\label{proposition:WMEP}
Consider $\mathcal X\subseteq\mathbb R^d$ and let $\P{X}{}\in\mathcal P_p(\mathcal X)$. Then the optimization problem in \eqref{equation:wasserstein_constrained_entropy} admits the dual formulation
\begin{customopti*}|s|
    {inf}
    {\alpha \geq 0, \varphi}
    {\log \int_\mathcal{X} \e{-\varphi(v)} \d{v} - \E{X}{\varphi^{c_\alpha}(x)} + \alpha \zeta^p}
    {}
    {}.
\end{customopti*}
The infimum is taken over $\alpha\ge 0$ and over admissible measurable dual potentials $\varphi:\mathcal X\to\mathbb R$ such that $\exp(-\varphi)$ is integrable and the terms in the dual objective are finite. The potential $\varphi$ is the Lagrange multiplier associated with the marginal consistency constraint between the entropy kernel and the transport kernel. Furthermore,
\[
c_\alpha(x,v):=\alpha\|x-v\|_p^p, \qquad \varphi^{c_\alpha}(x) :=
\Opt{inf}{v\in\mathcal{X}}{c_\alpha(x,v)-\varphi(v)}{}
\]
denotes the $c_\alpha$-transform of $\varphi$.

Moreover, suppose that strong duality holds and that an optimal dual pair
$(\alpha^W,\varphi^W)$ exists. Then the optimal induced marginal has
density
\[
p_V^W(v) = \frac{\e{-\varphi^W(v)}}{\int_{\mathcal X}\e{-\varphi^W(u)}\d{u}}.
\]
The corresponding entropy-maximizing conditional kernel is independent of $X$,
\[
\frac{\dP{V|X}{}^W}{\dP{V}{}^W} = 1, \qquad \P{X}{}\text{-a.s.}.
\]
\end{proposition}
The proof, including the derivation of the dual formulation and sufficient conditions for strong duality, is given in \cref{proof:proposition:WMEP}.

\begin{remark}[Entropy kernel versus transport coupling]
\label{remark:W-MEP}
At the W-MEP optimum the entropy-bearing kernel may be chosen independently of $X$, $\P{V|X=x}{}^{W}=\P{V}{}^{W}$ $\P{X}{}$-a.s. Therefore the entropy-bearing joint law can be factorized as $\P{X,V}{}^{W} = \P{X}{}\otimes\P{V}{}^{W}$, and $ \MutualInf{\P{X,V}{}^{W}}{X}{V} =0$. This does not imply that the optimal transport coupling between $\P{X}{}$ and $\P{V}{}$ is the product coupling. The former carries entropy, whereas the latter certifies Wasserstein feasibility.
\end{remark}

\begin{corollary}[W1-MEP]
\label{corollary:WMEP1}
Consider $\mathcal{X} = \mathbb{R}^d$, then, given $\P{X}{}  \in \mathcal{P}_1(\mathcal{X})$, the optimization problem in \eqref{equation:wasserstein_constrained_entropy}, with $p=1$, admits the dual formulation
\begin{customopti*}|s|
    {inf}
    {\varphi}
    {\log \int_{\mathbb{R}^d} \e{-\varphi(v)} \d{v} + \E{X}{\varphi(x)} + \emph{Lip}(\varphi) \zeta}
    {}
    {},
\end{customopti*}
where the infimum is restricted to pairs $(\alpha,\varphi)$ for which the log-partition term is finite, and $\text{Lip}_\alpha(\mathcal X)$ denotes the set of measurable $\alpha$-Lipschitz functions with respect to the $\ell_1$ norm, 
\[
\emph{Lip}_\alpha(\mathcal X) := \big\{ \varphi:\mathcal X\to\mathbb R: |\varphi(v)-\varphi(x)| \le \alpha\|v-x\|_1,\ \forall x,v\in\mathcal X \big\}.
\]
Moreover, if an optimal pair $(\alpha^W,\varphi^W)$ exists, then the optimal induced marginal has density
\[
p_V^W(v) = \frac{\e{-\varphi^W(v)}} {\int_{\mathcal X}\e{-\varphi^W(u)}\d{u}},
\]
and the entropy-maximizing conditional kernel is independent of $X$, $\P{V|X}{}^W=\P{V}{}^W, \ \P{X}{}\text{-a.s.}$
\end{corollary}

Proof provided in \cref{proof:corollary:WMEP1}.

\begin{remark}[$W_2$-MEP]
Consider $\mathcal X=\mathbb R^d$, and let $\P{X}{}\in\mathcal P_2(\mathcal X)$. Then the Wasserstein-constrained maximum entropy problem in \eqref{equation:wasserstein_constrained_entropy}, with $p=2$, admits the dual characterization
\begin{customopti*}|s|
    {inf}
    {\alpha \geq 0, \varphi}
    {\log \int_{\mathbb{R}^d} \e{-\varphi(v)} \d{v} + \E{X}{\varphi^{c}(x)} + \alpha \zeta^2}
    {}
    {},
\end{customopti*}
where the quadratic $c$-transform is defined by
\[
\varphi^c(x):= \Opt{inf}{v \in \mathcal{X}}{\alpha \|x-v\|_2^2-\varphi(v)}{}
\]
Moreover, suppose that strong duality holds and that an optimal dual pair $(\alpha^W,\varphi^W)$ exists. Then the corresponding entropy-maximizing conditional kernel is independent of $X$, and the optimal induced marginal has density as given in \cref{proposition:WMEP}.

This is the quadratic Wasserstein analogue of the Lipschitz-potential formulation in \cref{corollary:WMEP1} and is consistent with the JKO scheme \cite{JKO1998}, a variational Wasserstein formulation for the Fokker--Planck equation with decay. Furthermore, the W-MEP dual involves the quadratic $c_\alpha$-transform rather than a Lipschitz constraint. Under suitable regularity assumptions, the associated optimality conditions take the form of a Monge--Amp\`ere equation \cite{Villani2008,Santambrogio2015}. This contrasts with the $p=1$ case, where the dual reduces to a Lipschitz-potential formulation.
\end{remark}

\section{Cost of Dependence}
\label{section:gap}
The previous section characterized the optimal solutions of A-MEP and W-MEP, highlighting the different roles played by the perturbation kernel and the induced marginal distribution. We now investigate how these differences affect their optimal values. In particular, we seek to quantify the loss of optimality incurred when using the tractable A-MEP formulation as a surrogate for W-MEP. To this end, we first establish that the average-transport and Wasserstein constraints generate the same set of attainable marginal distributions. This allows us to express both problems over a common feasible set and identify the dependence between the nominal and perturbed variables, measured through mutual information, as the source of their discrepancy. We then use this characterization to derive bounds on the resulting entropy gap.

\begin{proposition}[Marginal equivalence]
\label{proposition:marginal_equivalence}
Let $\P{X}{}\in\mathcal{P}_p(\mathcal{X})$ be given, $\zeta\geq 0$, and let $F:\mathcal{P}_p(\mathcal{X})\to \mathbb{R}$ be any functional of the perturbed marginal distribution. Then,

\begin{minipage}{0.4\textwidth}
\begin{customopti*}|s|
    {sup}
    {\P{V|X}{} \in \mathcal{K}_p(\mathcal{X})}
    {F(\P{V}{})}
    {}
    {}
    \addConstraint{\E{X,V}{\|x-v\|_p^p} \leq \zeta^p}
\end{customopti*}
\end{minipage}
\begin{minipage}{0.05\textwidth}
$=$
\end{minipage}
\begin{minipage}{0.4\textwidth}
\begin{customopti*}|s|
    {sup}
    {\P{V}{} \in \mathcal{P}_p(\mathcal{X})}
    {F(\P{V}{}),}
    {}
    {}
    \addConstraint{W_p(\P{X}{}, \P{V}{}) \leq \zeta}.
\end{customopti*} 
\end{minipage}
\end{proposition}

\begin{proof}
`` $\leq$ '':
Let $\P{V|X}{}$ be a feasible argument for the problem on the left, with induced marginal $\P{V}{}$ and joint law $\P{X,V}{}$. Since $\P{X,V}{} \in \Pi(\P{X}{}, \P{V}{})$
\[
    W_p^p(\P{X}{},\P{V}{}) \leq \int_{\mathcal{X} \times \mathcal{X}} \|x-v\|_p^p \dP{X,V}{}(x,v) \leq\zeta^p.
\]
Hence, $\P{V}{}$ is feasible for the problem on the right and attains there the same objective value $F(\P{V}{})$. Since this holds for every feasible kernel on the left,

\begin{minipage}{0.4\textwidth}
\begin{customopti*}|s|
    {sup}
    {\P{V|X}{} \in \mathcal{K}_p(\mathcal{X})}
    {F(\P{V}{})}
    {}
    {}
    \addConstraint{\E{X,V}{\|x-v\|_p^p} \leq \zeta^p}
\end{customopti*}
\end{minipage}
\begin{minipage}{0.05\textwidth}
$\leq$
\end{minipage}
\begin{minipage}{0.4\textwidth}
\begin{customopti*}|s|
    {sup}
    {\P{V}{} \in \mathcal{P}_p(\mathcal{X})}
    {F(\P{V}{}),}
    {}
    {}
    \addConstraint{W_p(\P{X}{}, \P{V}{}) \leq \zeta}.
\end{customopti*} 
\end{minipage}

`` $\geq$ '': Conversely, let $\P{V}{}$ a feasible solution for the problem on the right, and let $\pi^\star\in\Pi(\P{X}{},\P{V}{}^\star)$ be an optimal coupling. Since $\mathcal{X}$ is Polish, $\pi^\star$ admits a disintegration
\[
    \dPQ{}{}^\star(x,v) =  \dP{X}{x} \dQ{V|X=x}{v}.
\]
Taking $\P{V|X}{} :=\Q{V|X}{}$, the resulting kernel induces the same marginal $\P{V}{}$ and satisfies
\[
    \int_{\mathcal{X}\times\mathcal{X}} \|x-v\|_p^p \d{(\P{X}{} \otimes \P{V|X}{})}(x,v) = W_p^p(\P{X}{},\P{V}{}) \leq\zeta^p.
\]
Hence, $\P{V|X}{}$ is feasible for the problem on the left and attains there the same objective value $F(\P{V}{})$,

\begin{minipage}{0.4\textwidth}
\begin{customopti*}|s|
    {sup}
    {\P{V|X}{} \in \mathcal{K}_p(\mathcal{X})}
    {F(\P{V}{})}
    {}
    {}
    \addConstraint{\E{X,V}{\|x-v\|_p^p} \leq \zeta^p}
\end{customopti*}
\end{minipage}
\begin{minipage}{0.05\textwidth}
$\geq$
\end{minipage}
\begin{minipage}{0.4\textwidth}
\begin{customopti*}|s|
    {sup}
    {\P{V}{} \in \mathcal{P}_p(\mathcal{X})}
    {F(\P{V}{}),}
    {}
    {}
    \addConstraint{W_p(\P{X}{}, \P{V}{}) \leq \zeta}.
\end{customopti*} 
\end{minipage}

\end{proof}

\cref{proposition:marginal_equivalence} applies to objectives that depend on the perturbation only through its induced marginal distribution. In particular, the two transport constraints are indistinguishable for any objective of the form $F(\P{V}{})$. When applying this result to entropy, we write $\ShannonEnt{}{\P{V}{}}$ and $\DiffEnt{}{\P{V}{}}$ for the Shannon and differential entropy functionals of the marginal measure, respectively. In particular,
\[
\ShannonEnt{}{\P{V}{}}=\ShannonEnt{\P{X,V}{}}{V},
\qquad
\DiffEnt{}{\P{V}{}}=\DiffEnt{\P{X,V}{}}{V}.
\]

\begin{corollary}[Common average-transport representation]
\label{corrollary:common_representation}
Let $A_\zeta$ and $W_\zeta$ denote the optimal values of A-MEP \ref{equation:average_constrained_entropy}, and W-MEP \ref{equation:wasserstein_constrained_entropy}, respectively. Then
\begin{customopti*}|s|
    {sup}
    {\P{V|X}{} \in \mathcal{K}_p(\mathcal{X})}
    {\DiffEnt{\P{X,V}{}}{V}}
    {}
    {W_\zeta =}
    {}
    \addConstraint{\E{X,V}{\|x-v\|_p^p} \leq \zeta^p}
\end{customopti*} 
whereas
\begin{customopti*}|s|
    {sup}
    {\P{V|X}{} \in \mathcal{K}_p(\mathcal{X})}
    {\DiffEnt{\P{X,V}{}}{V} - \MutualInf{\P{X,V}{}}{X}{V}}
    {}
    {A_\zeta = }
    \addConstraint{\E{X,V}{\|x-v\|_p^p} \leq \zeta^p,}
\end{customopti*}
and $A_\zeta\leq W_\zeta$.
\end{corollary}

\begin{proof}
The characterization of W-MEP established above can be equivalently written as 
\begin{customopti*}|s|
    {sup}
    {\P{V}{} \in \mathcal{P}_p(\mathcal{X})}
    {\DiffEnt{\P{X,V}{}}{V}}
    {}
    {W_\zeta =}
    {}
    \addConstraint{W_p(\P{X}{}, \P{V}{}) \leq \zeta.}
\end{customopti*} 
This is due to \cref{proposition:WMEP,remark:W-MEP} which show that the optimal argument for the W-MEP problem is an argument independent of $X$. Additionally, by \cref{proposition:marginal_equivalence} with $F(\P{V}{})=\DiffEnt{}{\P{V}{}}$, and using $\DiffEnt{}{\P{V}{}}=\DiffEnt{\P{X,V}{}}{V}$ results in,
\begin{customopti*}|s|
    {sup}
    {\P{V|X}{} \in \mathcal{K}_p(\mathcal{X})}
    {\DiffEnt{\P{X,V}{}}{V}}
    {}
    {W_\zeta =}
    {}
    \addConstraint{\E{X,V}{\|x-v\|_p^p} \leq \zeta^p}.
\end{customopti*} 
The expression for $A_\zeta$ follows from the identity
\[
    \DiffEnt{\P{X,V}{}}{V|X}=\DiffEnt{\P{X,V}{}}{V}-\MutualInf{\P{X,V}{}}{X}{V}.
\]
Finally, for every feasible kernel $\P{V|X}{}$, 
\[
\DiffEnt{\P{X,V}{}}{V} - \MutualInf{\P{X,V}{}}{X}{V} \le \DiffEnt{\P{X,V}{}}{V},
\]
by nonnegativity of mutual information. Suppose, by contradiction, that $A_\zeta>W_\zeta$. By the definition of supremum, there then exists a feasible kernel for which
\[
\DiffEnt{\P{X,V}{}}{V} - \MutualInf{\P{X,V}{}}{X}{V} >W_\zeta.
\]
For the same kernel, nonnegativity of mutual information implies
\[
\DiffEnt{\P{X,V}{}}{V} > W_\zeta,
\]
contradicting the definition of $W_\zeta$ as the supremum of $\DiffEnt{\P{X,V}{}}{V}$ over the common feasible set. Hence $A_\zeta\leq W_\zeta$.
\end{proof}

After expressing both problems over the same average-transport feasible set, their difference lies entirely in the objective. W-MEP maximizes the entropy of the induced marginal, whereas A-MEP maximizes the same marginal entropy net of the dependence $\MutualInf{\P{X,V}{}}{X}{V}$ required by the perturbation. We can now isolate the effect of the mutual-information term and characterize the resulting optimal-value gap, $\Delta := W_\zeta-A_\zeta \geq 0$.

\begin{theorem}[Information-gap representation and bounds]
\label{theorem:information_gap}
Suppose that A-MEP and W-MEP attain their optimal values. Let $\P{V|X}{}^A$ be an optimal kernel for A-MEP, with induced marginal $\P{V}{}^A$ and joint law $\P{X,V}{}^A = \P{X}{}\otimes\P{V|X}{}^A$, and let $\P{V}{}^W$ denote an optimal marginal associated with W-MEP. Define $I_A := \MutualInf{\P{X,V}{}^A}{X}{V}$, and 
\begin{customopti*}|s|
    {inf}
    {\P{V|X}{} \in \mathcal{K}_p(\mathcal{X})}
    {\MutualInf{\P{X,V}{}}{X}{V}}
    {}
    {I_W :=}
    {}
    \addConstraint{\Eonly{X} \E{V|X}{\|x-v\|_p^p} \leq \zeta^p}
    \addConstraint{\DiffEnt{\P{X,V}{}}{V} = W_\zeta}.
\end{customopti*} 
Then, the entropy gap $\Delta = \left[ \DiffEnt{\P{X,V}{}^W}{V}-\DiffEnt{\P{X,V}{}^A}{V} \right] +I_A$, is bounded by,
\[
    I_A \leq \Delta \leq I_W.
\]
\end{theorem}

\begin{proof}
Since $\P{V|X}{}^A$ is optimal for A-MEP, $A_\zeta = \DiffEnt{\P{X,V}{}^A}{V}-I_A$, while optimality of $\P{V}{}^W$ for W-MEP gives $W_\zeta=\DiffEnt{\P{X,V}{}^W}{V}$. Hence
\[
    \Delta = \DiffEnt{\P{X,V}{}^W}{V}-\DiffEnt{\P{X,V}{}^A}{V}+I_A.
\]
Since $\P{V}{}^{A}$ is attainable under the common average-transport constraint and $\P{V}{}^{W}$ maximizes marginal entropy over all attainable marginals,
\[
    \DiffEnt{\P{X,V}{}^W}{V} \geq \DiffEnt{\P{X,V}{}^A}{V},
\]
which proves $\Delta\geq I_A$.

For the upper bound, let $\P{V|X}{} \in \mathcal{K}_p$, with $\Eonly{X} \E{V|X}{\|x-v\|_p^p} \leq \zeta^p$, be any feasible kernel inducing the W-MEP marginal $\P{V}{}^W$. Since A-MEP maximizes $\DiffEnt{\P{X,V}{}}{V} - \MutualInf{\P{X,V}{}}{X}{V}$, 
\[
    A_\zeta \geq \DiffEnt{\P{X,V}{}}{\P{V}{}^W}-\MutualInf{\P{X,V}{}}{X}{V} =  W_\zeta-\MutualInf{\P{X,V}{}}{X}{V}.
\]
Therefore,
\[
  \DiffEnt{\P{X,V}{}^A}{V} - \DiffEnt{\P{X,V}{}^W}{V} \geq I_A -\MutualInf{\P{X,V}{}}{X}{V},
\]
and hence,
\[
    \Delta \leq \MutualInf{\P{X,V}{}}{X}{V}.
\]
Taking the infimum over all $\P{V|X}{}$ with $\Eonly{X} \E{V|X}{\|x-v\|_p^p} \leq \zeta^p$ inducing $\P{V}{}^W$ yields $\Delta \leq I_W$.
\end{proof}

\begin{corollary}[Explicit entropy-gap bounds on $\mathbb{R}^d$]
\label{corollary:euclidean-gap-bounds}
Let $\mathcal{X}=\mathbb{R}^d$, and suppose that the assumptions of \cref{theorem:information_gap} hold. Let $Z$ be independent of $X$ and distributed according to the generalized Gaussian identified in \cref{proposition:AMEP}, so that $V=X+Z$ is an optimal A-MEP kernel. Define
\[
m_p(X) := \left( \E{}{\|X-\E{}{X}\|_p^p} \right)^{1/p}, \qquad
B_\zeta(X) := d\log\left( 1+\frac{m_p(X)}{\zeta} \right).
\]
Then
\[
\MutualInf{\P{X,V}{}}{X}{X+Z} \leq \Delta \leq \min\left\{ I_W, B_\zeta(X) \right\} \leq B_\zeta(X).
\]
\end{corollary}

\begin{proof}
By \cref{proposition:AMEP}, the A-MEP optimum is attained by a generalized-Gaussian kernel centered at $X$, whose displacement distribution is independent of the nominal realization. The mutual information associated with this A-MEP optimizer is $I_A=\MutualInf{\P{X,V}{}}{X}{X + Z}$. Applying \cref{theorem:information_gap} gives
\[
\MutualInf{\P{X,V}{}}{X}{X+Z} \leq \Delta \leq I_W.
\]
It remains to establish the explicit upper bound. For any kernel $\P{V|X}{}$ feasible for W-MEP, by Minkowski's inequality we have, $\E{}{\|V-\E{}{X}\|_p^p} \leq (\zeta+m_p)^p$, since
\begin{align*}
\left( \E{}{\|V-\E{}{X}\|_p^p}\right) ^{1/p} &= \left( \E{}{\|(V-X)+(X-\E{}{X})\|_p^p} \right)^{1/p} \\
&\leq \left( \E{}{\|V-X\|_p^p} \right)^{1/p} + \left( \E{}{\|X-\E{}{X}\|_p^p} \right)^{1/p} \\
&\leq \zeta+ m_p.
\end{align*}

By the maximum-entropy characterization in \cref{proposition:AMEP}, and by translation invariance of differential entropy,
\[
\DiffEnt{\P{V}{}}{V} = \DiffEnt{\P{V}{}}{V-\E{}{X}} \leq A_{\zeta+m_p}.
\]
Since this holds for every feasible W-MEP kernel, then, $W_\zeta \leq A_{\zeta+m_p}$, and consequently,
\[
\Delta = W_\zeta-A_\zeta \leq A_{\zeta+m_p}-A_\zeta.
\]

Furthermore, from \cref{proposition:AMEP}, the A-MEP value scales with the distortion radius as
\[
A_r=C_{p,d}+d\log r,
\]
where $C_{p,d}$ does not depend on $r$. Hence
\[
\Delta \leq d\log\left( \frac{\zeta+m_p}{\zeta} \right) = d\log\left( 1+\frac{m_p}{\zeta} \right) 
= B_\zeta(X).
\]
Combining this bound with the result from \cref{theorem:information_gap}, $\Delta\leq I_W$, gives
\[
\Delta \leq \min\left\{ I_W, B_\zeta(X) \right\} \leq B_\zeta(X),
\]
which completes the proof.
\end{proof}

\begin{remark}[Dependence of the gap bound on the transport radius]
The explicit upper bound $B_\zeta(X) = d\log\left(1+\frac{m_p(X)}{\zeta}\right)$ is decreasing in the transport radius $\zeta$. In particular, $B_\zeta(X)\to 0 \quad\text{as}\quad \zeta\to\infty$,  whereas $B_\zeta(X)\to\infty$ as $\zeta\downarrow 0$ whenever $m_p(X)>0$. The monotonicity of $B_\zeta(X)$ does not imply that the entropy gap $\Delta_\zeta$ is itself monotone in $\zeta$, since an ordering of the upper bounds does not induce an ordering of the corresponding gaps.

This tells that for small transport budgets, A-MEP requires the entropy-bearing perturbation to remain close to each realization of $X$, thereby retaining dependence between $X$ and $V$. In contrast, W-MEP constrains only the induced marginal distribution and allows the entropy-bearing kernel to be independent of $X$. As the transport budget increases, this distinction becomes
less restrictive. In particular, \cref{corollary:euclidean-gap-bounds} implies that $\Delta_\zeta\to 0 \quad\text{as}\quad \zeta\to\infty$, provided that $m_p(X)<\infty$.
\end{remark}
\section{Discretely Supported Nominal Distribution}
\label{section:discrete}

A case of particular interest is that of a discretely supported nominal distribution. Such distributions arise naturally in data-driven settings, where the underlying distribution is unknown and is represented by an empirical distribution constructed from a finite collection of observations. In this case, the nominal law takes the form
\begin{equation}
\label{equation:discrete_nominal}
\P{X}{} = \sum_{i=1}^n a_i\delta_{x_i}, \qquad a_i>0,\qquad \sum_{i=1}^n a_i=1,
\end{equation}
with the empirical distribution corresponding to the particular choice $a_i=1/n$. Although the nominal distribution is discretely supported, the perturbed variable $V$ is allowed to take values in the ambient space $\mathbb{R}^d$. The resulting formulation is more general than a pointwise perturbation of the nominal support since the probability mass associated with each $x_i$ may be split and transported according to an arbitrary distribution on $\mathbb{R}^d$. In particular, translations of the nominal support points are recovered as the special case in which the corresponding conditional distributions are Dirac measures. In this setting, the A-MEP and W-MEP problems admit particularly transparent representations, highlighting the distinction between the entropy of the conditional perturbations and the entropy of their induced mixture.

\begin{proposition}[Semi--discrete representation]
\label{proposition:semi_discrete}
Let $\mathcal{X}=\mathbb{R}^d$ and consider an empirical, finitely supported nominal distribution as in \eqref{equation:discrete_nominal}. Then the A-MEP problem \ref{equation:average_constrained_entropy} is equivalent to
\begin{customopti*}|s|
    {sup}
    {\{\P{i}{} \in \mathcal{P}_p(\mathbb{R}^d) \}_{i=1}^n}
    {\sum_{i=1}^n a_i \DiffEnt{\P{i}{}}{V}}
    {}
    {A = }
    \addConstraint{\sum_{i=1}^n a_i \int_{\mathbb{R}^d} \|x_i-v\|_p^p \dP{i}{v} \leq \zeta^p,}
\end{customopti*}
and  the W-MEP problem \ref{equation:wasserstein_constrained_entropy} admits the equivalent average-transport representation
\begin{customopti*}|s|
    {sup}
    {\{\P{i}{} \in \mathcal{P}_p(\mathbb{R}^d) \}_{i=1}^n}
    {\DiffEnt{\overline{\P{}{}}}{V}}
    {}
    {W = }
    \addConstraint{\sum_{i=1}^n a_i \int_{\mathbb{R}^d} \|x_i-v\|_p^p \dP{i}{v} \leq \zeta^p,}
\end{customopti*}
where $\P{i}{}=\P{V|X=x_i}{},\ i=1,\ldots,n$ and, $\overline{\P{}{}}:=\sum_{i=1}^n a_i\P{i}{}$ denotes the induced $V$-marginal.
Moreover, the quantities $\DiffEnt{\P{i}{}}{V}$ and $\DiffEnt{\overline{\P{}{}}}{V}$  denote the differential entropy of $V$ under $\P{i}{}$ and $\overline{\P{}{}}$, respectively, whenever well defined.
\end{proposition}

\begin{proof}
Let the nominal distribution be discretely supported on as given by \eqref{equation:discrete_nominal}. Although $\P{X}{}$ is discretely supported, the perturbed variable $V$ is allowed to take values in the ambient space $\mathbb{R}^d$. Therefore, for each $x_i$, define
\[
\P{i}{} := \P{V|X=x_i}{} \in\mathcal{P}_p(\mathbb{R}^d), \qquad i=1,\ldots,n.
\]
The collection $\{\P{i}{}\}_{i=1}^n$ completely determines the perturbation kernel on the support of $\P{X}{}$. The joint distribution induced by this kernel is
\[
\P{X,V}{x,v} = \sum_{i=1}^n a_i \delta_{x_i}(x) \dP{i}{v},
\]
and its $V$-marginal is therefore the mixture $\overline{\P{}{}} = \sum_{i=1}^n a_i \P{i}{}$. Since $X$ is discrete and $V$ is continuous, the conditional differential entropy is obtained by averaging the differential entropies of the conditional laws,
\[
\DiffEnt{\P{X,V}{}}{V|X} = \sum_{i=1}^n \P{X}{x_i} \DiffEnt{\P{i}{}}{v} = \sum_{i=1}^n a_i \DiffEnt{\P{i}{}}{V}.
\]
If $\P{i}{}$ admits a density $b_i$ with respect to Lebesgue measure, this can equivalently be written as
\[
\DiffEnt{\P{X,V}{}}{V|X} = - \sum_{i=1}^n a_i \int_{\mathbb{R}^d} b_i(v)\log b_i(v) \d{v}.
\]
Furthermore, the average transport cost under the perturbation-induced joint law is
\[
\Eonly{X}\E{V|X}{\|X-V\|_p^p} = \sum_{i=1}^n a_i \int_{\mathbb{R}^d} \|x_i-v\|_p^p \dP{i}{v}.
\]
As a result, the A-MEP problem is equivalent to
\begin{customopti*}|s|
    {sup}
    {\{\P{i}{} \in\mathcal{P}_p(\mathbb{R}^d)\}_{i=1}^n}
    {\sum_{i=1}^n a_i \DiffEnt{\P{i}{}}{V}}
    {}
    {}
    \addConstraint{\sum_{i=1}^n a_i \int_{\mathbb{R}^d} \|x_i-v\|_p^p \dP{i}{v} \leq \zeta^p
    }.
\end{customopti*}

For W-MEP, let us recall the equivalent average-transport representation of \cref{corrollary:common_representation},
\begin{customopti*}|s|
    {sup}
    {\P{V|X}{} \in \mathcal{K}_p(\mathcal{X})}
    {\DiffEnt{\P{X,V}{}}{V}}
    {}
    {}
    \addConstraint{\Eonly{X} \E{V|X}{\|x-v\|_p^p} \leq \zeta^p},
\end{customopti*} 
Since $V$-marginal is the mixture $\overline{\P{}{}} = \sum_{i=1}^n a_i \P{i}{}$, the W-MEP result follows directly,
\begin{customopti*}|s|
    {sup}
    {\{\P{i}{} \in\mathcal{P}_p(\mathbb{R}^d)\}_{i=1}^n}
    {\DiffEnt{\overline{\P{}{}}}{V}}
    {}
    {}
    \addConstraint{\sum_{i=1}^n a_i \int_{\mathbb{R}^d} \|x_i-v\|_p^p \dP{i}{v} \leq \zeta^p
    }.
\end{customopti*}

\end{proof}

Under a discretely supported nominal distribution, W-MEP maximizes the differential entropy of the mixture induced by the conditional transport distributions $\{\P{i}{}\}_{i=1}^n$, subject to the same average transport constraint as A-MEP. In contrast, A-MEP maximizes the average of their individual differential entropies, whereas W-MEP maximizes the entropy of the average marginals.

\begin{corollary}[Fully discrete representation]
\label{corollary:fully_discrete}
Suppose that the perturbed variable is restricted to the finite set $\{v_1,\ldots,v_m\}$ and consider an empirical, finitely supported nominal distribution as in \eqref{equation:discrete_nominal}. Define, for each $i=1,\ldots,n$, 
\[
    b_{ij} := \P{V\mid X=x_i}{v_j}, \qquad j=1,\ldots,m,
\]
and let $b_i=(b_{i1},\ldots,b_{im}) \in \Delta_m$ where $\Delta_m = \left\{ b\in\mathbb{R}_+^m: \sum_{j=1}^m b_j=1 \right\}$. Then the A-MEP problem admits the finite-dimensional representation
\begin{customopti*}|s|
    {sup}
    {\{b_i \in \Delta_m\}_{i=1}^n}
    {\sum_{i=1}^n a_i \ShannonEnt{}{b_i}}
    {}
    {}
    \addConstraint{\sum_{i=1}^n \sum_{j=1}^m a_i b_{ij} \|x_i-v_j\|_p^p \leq \zeta^p},
\end{customopti*}
with $\ShannonEnt{}{b_i} = -\sum_{j=1}^m b_{ij}\log b_{ij}$. Moreover, the W-MEP problem admits the equivalent representation
\begin{customopti*}|s|
    {sup}
    {\{b_i \in \Delta_m\}_{i=1}^n}
    {\ShannonEnt{}{\sum_{i=1}^n a_i b_i}}
    {}
    {}
    \addConstraint{\sum_{i=1}^n \sum_{j=1}^m a_i b_{ij} \|x_i-v_j\|_p^p \leq \zeta^p}.
\end{customopti*} 
\end{corollary}

\begin{proof}
This is the fully discrete specialization of \cref{proposition:semi_discrete}. Indeed, since $V\in\mathcal{V}$, each conditional law $b_i=\P{V\mid X=x_i}{}$ is identified with a probability vector $b_i\in\Delta_m$, and its entropy is $\ShannonEnt{}{b_i} = -\sum_{j=1}^m b_{ij}\log b_{ij}$. The induced marginal satisfies $\P{V}{v_j} = \sum_{i=1}^n a_i b_{ij}$, and the average transport cost becomes $\sum_{i=1}^n a_i \sum_{j=1}^m b_{ij}\|x_i-v_j\|_p^p$. Substituting these expressions into the A-MEP and W-MEP representations of \cref{proposition:semi_discrete}
gives the fully discrete representations of A-MEP and W-MEP.
\end{proof}

\subsection{Entropy gap}

Consider the empirical, finitely supported nominal distribution, in \eqref{equation:discrete_nominal}, and let $\{\P{i}{}^A\}_{i=1}^n$ and $\{\P{i}{}^W\}_{i=1}^n$ denote optimal conditional perturbation laws for A-MEP and W-MEP, respectively. Define the corresponding induced marginals,
\[
\overline{\P{}{}}^A:= \sum_{i=1}^n a_i \P{i}{}^A, \qquad \overline{\P{}{}}^W := \sum_{i=1}^n a_i \P{i}{}^W,
\]
and suppose that the relevant differential entropies are finite. The following proposition characterizes the gap between the results in \cref{proposition:semi_discrete}, 
\[
\Delta = \DiffEnt{\overline{\P{}{}}^W }{V} - \sum_{i=1}^n a_i \DiffEnt{\P{i}{}^A}{V}.
\]

\begin{proposition}[Decomposition of the entropy gap]
\label{proposition:gap_decomposition}
The entropy gap $\Delta$ admits an exact decomposition,
\[
\Delta = \DiffEnt{\overline{\P{}{}}^W}{V}-\DiffEnt{\overline{\P{}{}}^A}{V} + \sum_{i=1}^n a_i \KL{\P{i}{}^A}{\overline{\P{}{}}^A},
\]
equivalently given by,
\[
\Delta = \DiffEnt{\overline{\P{}{}}^W}{V}-\DiffEnt{\overline{\P{}{}}^A}{V} + \MutualInf{\P{X,V}{}^A}{X}{V}.
\]
\end{proposition}

\begin{proof}
From the finite-support representations of A-MEP, $A = \sum_{i=1}^n a_i \DiffEnt{\P{i}{}^A}{V}$, and W-MEP, $W = \DiffEnt{\overline{\P{}{}}^W}{V}$. Hence $\Delta = W-A = \DiffEnt{\overline{\P{}{}}^W}{V} - \sum_{i=1}^n a_i \DiffEnt{\P{i}{}^A}{V}$. Adding and subtracting $\DiffEnt{\overline{\P{}{}}^A}{V}$ gives
\[
\Delta = \DiffEnt{\overline{\P{}{}}^W}{V}-\DiffEnt{\overline{\P{}{}}^A}{V} + \DiffEnt{\overline{\P{}{}}^A}{V} - \sum_{i=1}^n a_i \DiffEnt{\P{i}{}^A}{V}.
\]
It remains to identify the second term. Let $p_i^A$ denote the density
of $\P{i}{}^A$ with respect to Lebesgue measure and define $\overline{p}^A(v) := \sum_{i=1}^n a_i p_i^A(v)$, which is the density of $\overline{\P{}{}}^A$. Then
\[
\begin{aligned}
\sum_{i=1}^n a_i \KL{\P{i}{}^A}{\overline{\P{}{}}^A} &= \sum_{i=1}^n a_i  \int_{\mathbb R^d} p_i^A(v) \log \frac{p_i^A(v)}{\overline{p}^A(v)} \d{l} \\
&= \sum_{i=1}^n a_i \int_{\mathbb R^d} p_i^A(v)\log p_i^A(v) \d{v} \\
&\quad - \sum_{i=1}^n a_i \int_{\mathbb R^d} p_i^A(v)\log \overline{p}^A(v)\d{v}.
\end{aligned}
\]
Using $\sum_{i=1}^n a_i p_i^A(v)=\overline{p}^A(w)$, the second term becomes
\[
\sum_{i=1}^n a_i \int_{\mathbb R^d} p_i^A(v)\log \overline{p}^A(v) \d{v} = \int_{\mathbb R^d} \overline{p}^A(v)\log \overline{p}^A(v) \d{v}.
\]
Therefore,
\begin{align*}
\sum_{i=1}^n a_i \KL{\P{i}{}^A}{\overline{\P{}{}}^A} &= \sum_{i=1}^n a_i \int_{\mathbb R^d} p_i^A(v)\log p_i^A(v)\d{v} - \int_{\mathbb R^d} \overline{p}^A(v)\log \overline{p}^A(v)\d{v} \\
&= \DiffEnt{\overline{\P{}{}}^A}{V} - \sum_{i=1}^n a_i \DiffEnt{\P{i}{}^A}{V}.
\end{align*}
Substituting this identity into the expression for $\Delta$ yields
\[
\Delta = \DiffEnt{\overline{\P{}{}}^W}{V} - \DiffEnt{\overline{\P{}{}}^A}{V} + \sum_{i=1}^n a_i \KL{\P{i}{}^A}{\overline{\P{}{}}^A}.
\]
Finally, under the A-MEP joint law,
\[
\P{X,V}{}^A(x,v) = \sum_{i=1}^n a_i \delta_{x_i}(x) \dP{i}{}^A(v),
\]
the conditional law of $V$ given $X=x_i$ is $\P{i}{}^A$, while the $V$-marginal is $\overline{\P{}{}}^A$. Hence
\[
\MutualInf{\P{X,V}{}^A}{X}{V} = \sum_{i=1}^n a_i \KL{\P{i}{}^A}{\overline{\P{}{}}^A},
\]
which proves the equivalent mutual-information representation.
\end{proof}

\begin{remark}[Jensen--Shannon divergence]
The last term is the equal-weight Jensen--Shannon divergence \cite{Lin1991},
\[
\operatorname{JS}(\P{1}{}^A,\ldots,\P{n}{}^A) := \sum_{i=1}^n a_i\KL{\P{i}{}^A}{a_i \sum_{j=1}^n \P{j}{}^A}
\]
and therefore 
\[
\Delta = \DiffEnt{\overline{\P{}{}}^W}{V}-\DiffEnt{\overline{\P{}{}}^A}{V} + \operatorname{JS}(\P{1}{}^A,\ldots,\P{n}{}^A).
\] 
The decomposition separates two sources of discrepancy between A-MEP and W-MEP. The first term, $\DiffEnt{\overline{\P{}{}}^W}{V}-\DiffEnt{\overline{\P{}{}}^A}{V}$, measures the gain in marginal entropy obtained by optimizing directly over the perturbed marginal. The second term is the mutual information retained by the A-MEP perturbation mechanism, or equivalently the weighted Jensen--Shannon divergence among its conditional perturbation laws. 
\end{remark}

\section{A symmetric two-point example}
We conclude this section with a simple one-dimensional example illustrating the entropy-gap decomposition in \cref{proposition:gap_decomposition}. Let $p=2$, and consider, on $\mathbb{R}$, the symmetric nominal distribution
\[
    \P{X}{} =  \frac{1}{2}\delta_{-a} + \frac{1}{2}\delta_a.
\]
For A-MEP, \cref{corollary:AMEP2} gives the optimal additive Gaussian perturbation $V = X + Z$, $Z \sim \mathcal{N}(0,\zeta^2)$ and $Z \perp X$. Hence, the two optimal conditional laws are
\[
    \P{-}{}^A = \mathcal{N}(-a,\zeta^2), \qquad \P{+}{}^A = \mathcal{N}(a,\zeta^2),
\]
and the corresponding induced marginal is the Gaussian mixture
\[
    \P{V}{}^A = \frac{1}{2}\P{-}{}^A + \frac{1}{2}\P{+}{}^A = \frac{1}{2}\mathcal{N}(-a,\zeta^2) + \frac{1}{2}\mathcal{N}(a,\zeta^2).
\]
The A-MEP optimal value follows directly from \cref{corollary:AMEP2},
\[
    A_\zeta = \DiffEnt{\P{X,V}{}}{V|X} = \frac{1}{2} \log (2\pi e \zeta^2).
\]
In contrast, W-MEP maximizes the entropy of the induced marginal over the Wasserstein ball. Since
$\P{X}{}$ is symmetric around the origin, an optimal W-MEP marginal may be chosen symmetric. 
Indeed, for any feasible marginal $\P{V}{}$, let $\P{V_-}{}$ denote the distribution of $-V$ and let
\[
\Q{}{} =\frac{1}{2} \P{V}{} + \frac{1}{2} \P{V_-}{}.
\]
Since $X$ and $-X$ have the same distribution and the quadratic transport cost is invariant under reflection,
\[
W_2(\P{X}{},\P{V_-}{})=W_2(\P{X}{},\P{V}{}).
\]
Moreover, $W_2^2(\P{X}{},\cdot)$ is convex with respect to mixtures, so
\[
W_2^2(\P{X}{},\Q{}{}) \le \frac{1}{2} W_2^2(\P{X}{},\P{V}{}) + \frac{1}{2} W_2^2(\P{X}{},\P{V_-}{})
= W_2^2(\P{X}{},\P{V}{}).
\]
Hence, if $\P{V}{}$ is feasible, then $\Q{}{}$ is also feasible. On the other hand, differential entropy is concave and invariant under reflection, and therefore
\[
\DiffEnt{}{\Q{}{}} \ge \frac{1}{2} \DiffEnt{}{\P{V}{}} + \frac{1}{2} \DiffEnt{}{\P{V_-}{}} = \DiffEnt{}{\P{V}{}}.
\]
Thus, symmetrization does not increase the Wasserstein cost and does not decrease the entropy. Consequently, an optimal marginal may be taken symmetric about the origin. Furthermore, the Wasserstein distance can be characterized explicitly using the cumulative distribution function \cite[Proposition 2.17]{Santambrogio2015},
\[
W_2^2(\P{X}{},\P{V}{}) = \int_0^1 \left| F_X^{-1}(t)-F_V^{-1}(t) \right|^2 \d{t},
\]
where $F^{-1}$ denotes the generalized inverse of the corresponding cumulative distribution function. In the particular case in this section, the generalized inverse is
\[
F_X^{-1}(t) = 
\begin{cases}
-a, & 0<t<\frac{1}{2},\\
 a, & \frac{1}{2}<t<1.
\end{cases}
\]
Therefore,
\[
W_2^2(\P{X}{},\P{V}{}) = \int_0^{1/2} \left(F_V^{-1}(t)+a\right)^2 \d{t} + \int_{1/2}^1 \left(F_V^{-1}(t)-a\right)^2 \d{t}.
\]
Since $\P{V}{}$ may be taken symmetric around the origin, $F_V^{-1}(1-t)=-F_V^{-1}(t)$ for $t\in(0,1)$. Changing variables $t\mapsto 1-t$ in the first integral, 
\[
W_2^2(\P{X}{},\P{V}{}) = 2\int_{1/2}^1 \left(F_V^{-1}(t)-a\right)^2 \d{t} = \int_0^1 \left(\left|F_V^{-1}(t)\right|-a\right)^2 \d{t},
\]
where the equality holds by the symmetry argument. Alternativelly, in the usual measure notation, 
\[
W_2^2(\P{X}{},\P{V}{}) = \int_{\mathbb R} \left(|v|-a\right)^2 \dP{V}{v}.
\]
As a result, the W-MEP simplifies to,
\begin{customopti*}|s|
    {sup}
    {\P{V}{} \in \mathcal{P}_2(\mathcal{X})}
    {\DiffEnt{}{V}}
    {}
    {}
    \addConstraint{\E{V}{\left(|v|-a\right)^2} \leq \zeta^2.}
\end{customopti*} 

\begin{figure}[t!]
    \centering
    \begin{tikzpicture}
    \begin{axis}[
        width=0.82\textwidth,
        height=6.5cm,
        domain=-3.5:3.5,
        samples=500,
        xlabel={$v$},
        ylabel={Density},
        xmin=-3.5,
        xmax=3.5,
        ymin=0,
        ymax=0.4,
        axis lines=left,
        tick align=outside,
        legend style={
            anchor=south,
            legend columns=1,
            draw=none
        },
    ]


    \addplot[
        thick,
        blue
    ]
    {
        0.5/(sqrt(2*pi)*0.6)
        * exp(-(x+1)^2/(2*0.6^2))
        +
        0.5/(sqrt(2*pi)*0.6)
        * exp(-(x-1)^2/(2*0.6^2))
    };
    \addlegendentry{$A-MEP$}

    \addplot[
        thick,
        red
    ]
    {
        exp(-1.093892*(abs(x)-1)^2)/3.131
    };
    \addlegendentry{$W-MEP$}

    \addplot[
        dashed,
        thin,
        forget plot
    ] coordinates {(-1,0) (-1,0.72)};

    \addplot[
        dashed,
        thin,
        forget plot
    ] coordinates {(1,0) (1,0.72)};

    \node[below] at (axis cs:-1,0) {$-a$};
    \node[below] at (axis cs:1,0) {$a$};

    \end{axis}
    \end{tikzpicture}
    \caption{
Optimal induced marginals under A-MEP and W-MEP for $\P{X}{}=\frac{1}{2}\delta_{-a}+\frac{1}{2}\delta_a$, with $a=1$ and $\zeta=0.6$. The A-MEP marginal is the Gaussian mixture $\P{V}{}^A =\frac{1}{2}\mathcal{N}(-a,\zeta^2) +\frac{1}{2}\mathcal{N}(a,\zeta^2)$, whereas the W-MEP marginal has density $p_V^{W}(v)\propto \exp[-\alpha^\star(|v|-a)^2]$, with $\alpha^\star\approx1.094$ determined by the active transport constraint. Dashed lines indicate the nominal support points $\{-a,a\}$.
}
\label{figure:comparison}
\end{figure}
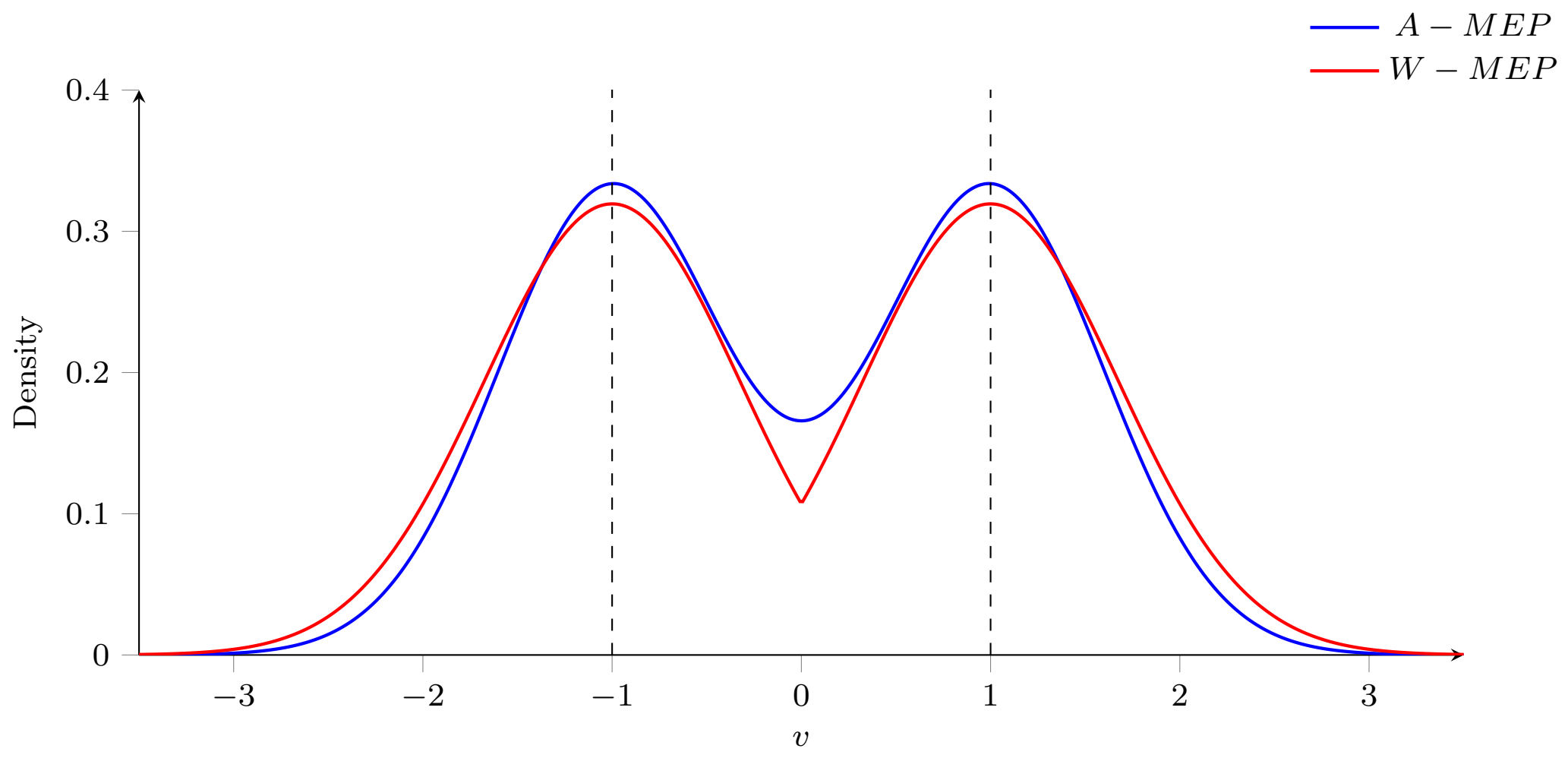

The problem is now a standard maximum-entropy problem under an expectation constraint. By the same Gibbs variational argument used in \cref{proposition:AMEP}, its optimal density is therefore of the form
\[
p_V^{W}(v) = \frac{1}{Z(\alpha)} \e{-\alpha (|v|-a)^2}, \qquad \alpha>0,
\]
with,
\[
Z(\alpha) = \int_{\mathbb{R}} \e{-\alpha (|v|-a)^2} \d{v} = \sqrt{\frac{\pi}{\alpha}} \left(1+\operatorname{erf}(a\sqrt{\alpha})\right).
\]
Furthermore, the optimal multiplier $\alpha^\star$ is chosen so that the transport constraint is active,
\[
\E{V}{(|v|-a)^2} = \zeta^2.
\]

\cref{figure:comparison} illustrates the resulting difference between the A-MEP and W-MEP optimal marginals for a representative choice of $a$ and $\zeta$.

The corresponding W-MEP optimal value is
\[
W_\zeta=\log Z(\alpha^\star)+\alpha^\star\zeta^2,
\]
and hence the entropy gap is
\[
\Delta_\zeta = \log Z(\alpha^\star)+\alpha^\star\zeta^2 -\frac12\log(2\pi e\zeta^2).
\]

This strict gap provides a concrete illustration of the distinction motivating the two formulations. A-MEP constrains the transition from each nominal realization $X$ to its perturbation $V$, and the resulting optimal kernel therefore retains information about the nominal support point from which $V$ originates. W-MEP, in contrast, constrains only the perturbed marginal and does not require the entropy-bearing mechanism to preserve this individual correspondence. The additional freedom at the distributional level allows W-MEP to attain a strictly larger entropy, as quantified here by $\Delta_\zeta>0$.

\section{Conclusion}

This paper investigated a tractable approximation to Wasserstein-constrained maximum-entropy inference through an average-transport formulation. We showed that A-MEP admits an explicit generalized-Gaussian additive-noise solution on $\mathbb{R}^d$, providing a computationally simple construction of a feasible distribution for W-MEP. We established quantitative bounds on the difference between their optimal values, including an explicit upper bound depending on the transport radius, dimension, and dispersion of the nominal distribution. This bound vanishes as the transport radius increases, providing a regime in which A-MEP approximates the optimal value of W-MEP arbitrarily closely.

The approximation guarantees follow from an information-theoretic characterization of the two formulations. We showed that A-MEP and W-MEP admit a common feasible-marginal representation, under which W-MEP maximizes marginal entropy, whereas A-MEP additionally penalizes the mutual information between the nominal and perturbed variables. We also derived a dual characterization of W-MEP in terms of Kantorovich potentials and, for finitely supported nominal distributions, an exact decomposition of the entropy gap into a marginal-entropy gain and a weighted Jensen--Shannon divergence. A symmetric two-point example illustrated that the gap can be strict. Overall, these results establish A-MEP as a tractable surrogate for W-MEP with quantitative optimality guarantees, while identifying the retained dependence between nominal and perturbed variables as the source of their discrepancy.

\bibliographystyle{styles/spmpsci}
\bibliography{main}

\appendix
\section{Proofs}
\subsection{Proof of AMEP results}

To prepare for the main result, we rewrite the conditional entropy in a form that separates its dependence on the marginal and conditional laws. This representation will serve as a convenient starting point for the analysis that follows.

\begin{lemma}[Equivalent representation of conditional entropy] 
\label{lemma:equivalent_formulation} 
For any joint distribution $\P{X,V}{}$ such that $\P{V}{}$ admits a density $p(v)$, the conditional entropy $\DiffEnt{\P{X,V}{}}{V|X}$ admits the representation,
\begin{equation*}
    \E{X}{\E{V|X}{-\log p(v)} - \KL{\P{V|X}{}}{\P{V}{}}}.
\end{equation*}
\end{lemma}

\begin{proof}
\label{proof:lemma:equivalent_formulation}
\begin{equation*}
\begin{aligned}
     \DiffEnt{\P{X,V}{}}{V|X} =& \int -\log p(v|x) p(v,x) \d{v} \d{x}\\
     =& \int -\log \frac{p(v,x)}{p(x)} p(v,x) \d{v} \d{x}\\
     =& \int \int \Big(- \log p(v) p(v|x) - \log \frac{p(v|x)}{p(v)} p(v|x) \Big) \d{v} \ p(x) \d{x},
\end{aligned}
\end{equation*}
and we have the result,
\begin{equation*}
     \E{X}{\E{V|X}{-\log p(v)} - \KL{\P{V|X}{}}{\P{V}{}}}.
\end{equation*}

Note that these also hold in a discrete case, since $\DiffEnt{\P{X,V}{}}{V|X} = \sum_{x,v} -\log p(v|x) p(v,x)$ and all $\log$ manipulations would equally hold.
\end{proof}

This representation prepares the ground for a variational argument based on the Donsker--Varadhan variational formula (also known as the Gibbs variational principle) \cite{Gibbs1902, DonskerVaradhan1975, DonskerVaradhan1983, DupuisEllis1997, DemboZeitouni1998, WainwrightJordan2008},  

\begin{lemma}[Donsker--Varadhan/Gibbs variational formula]
\label{lemma:DV}
Let $(\Theta,\mathcal{F})$ be a measurable space and let $\P{}{}$ be a probability
measure on it. For any measurable function $f:\Theta\to\mathbb{R}$ such that
$\E{\P{}{}}{e^{f}}<\infty$ and any probability measure $\Q{}{}$ with $\Q{}{} \ll \P{}{}$, the following dual variational representations hold:
\[
\log \E{\P{}{}}{e^{f}} = \Opt{sup}{\Q{}{}}{\E{\Q{}{}}{f - \KL{\Q{}{}}{\P{}{}}}}{} 
\qquad
\text{and}
\qquad
\KL{\Q{}{}}{\P{}{}} = \Opt{sup}{f}{\E{\Q{}{}}{f} - \log \E{\P{}{}}{e^{f}}}{}.
\]
Moreover, the supremum in the first identity is attained at the Gibbs measure
$\Q{}{}^\star$ defined by
\[
\frac{\dQ{}{}^\star}{\dP{}{}} (\theta) = \frac{e^{f(\theta)}}{\E{\vartheta \sim \P{}{}}{e^{f(\vartheta)}}}.
\]
and the supremum in the second identity is attained (up to an additive constant)
at
\[
f^\star(\theta)=\log\frac{\dQ{}{}}{\dP{}{}}(\theta).
\]
The two identities are convex conjugates of each other in the sense of
Fenchel--Legendre duality.
\end{lemma}
This identity provides a variational representation of log-partition functions in terms of KL divergence. It also yields exponentially tilted optimizers, which we exploit to derive the Gibbs structure under Wasserstein constraints.

The following proposition provides an exact formula for the integral $\int_{\mathbb{R}^d}\e{-\alpha\|x-v\|_p^p}\d{v}$ in terms of the unit $\ell_p$-ball volume and the Gamma function. This identity will be used in the proof of the main theorem to verify integrability and to simplify the resulting expressions.

\begin{proposition}[Integral for the $q$-radial kernel]
\label{proposition:lq_radial_integral}
Let $d\in\mathbb{N}$, $p \geq 1$, and $\alpha > 0$. For any $x\in\mathbb{R}^d$ the following integral is finite, independent of $x$, and admits the closed form
\[
\int_{\mathbb{R}^d}\e{-\alpha\|x-v\|_p^{p}} \d{v} = d \mathrm{vol}(B_p^d) \frac{\Gamma(d/p)}{p \alpha^{d/p}}
\]
where $\|\cdot\|_p$ denotes the $\ell_p$ norm and the integral is with respect to Lebesgue measure on $\mathbb{R}^d$, $\Gamma(\cdot)$ is the Gamma function, $\Gamma(t) = \int_{0}^{\infty} s^{t-1}e^{-s}\,ds$, $t>0$,
and $B_p^d=\{u\in\mathbb{R}^d:\|u\|_p\le 1\}$ is the $\ell_p$ unit ball with
$\mathrm{vol}(B_p^d)=\frac{\bigl(2\Gamma(1+\frac{1}{p})\bigr)^d}{\Gamma(1+\frac{d}{p})}$.
\end{proposition}

\begin{proof}
\label{proof:proposition:lq_radial_integral}
Note that the integral, 
\begin{equation*}
\int_{\mathbb{R}^d} \e{- \alpha \|x-v\|^p_p} \d{v},
\end{equation*}
takes the form
\begin{equation*}
    \int_{\mathbb{R}^d} \e{- \alpha \|u\|^p_p} \d{u}.
\end{equation*}
In the particular case of $p=2$, we can use \cite[Prop 2.51]{Book_Folland1999}, to write the above as,
\[
    \sigma(S^{d-1}) \int_0^\infty \e{- \alpha r^p} r^{d-1} \d{r}.
\]
That being said, this decomposition does not hold for a general $\ell_p$-ball. In this case we use \cite[Prop. 1]{Naor2003}, which tells us that, if $p \geq 1$, for $B_p^d:=\{u\in\mathbb{R}^d:\|u\|_p\le 1\}$ and $\partial B_p^d=\{u:\|u\|_p=1\}$, with $\mu_p$ as the \emph{cone probability measure} on $\partial B_p^d$, i.e.
$\mu_p(A)=\frac{\mathrm{vol}\big(\{ r z:\; 0\le r\le 1,\ z\in A\}\big)}{\mathrm{vol}(B_p^d)}$ for Borel sets $A\subseteq \partial B_p^d$, then for, every integrable function $f:\mathbb{R}^d\to\mathbb{R}$,
\begin{equation}
\label{equation:naor_decomp}
\int_{\mathbb{R}^d} f(u)\d{u} = d\,\mathrm{vol}(B_p^d)\int_0^\infty r^{d-1} \left(\int_{\partial B_p^d} f(r z)\d{\mu_p(z)}\right)\d{r},
\end{equation}
with the volume of the unit $\ell_p$-ball, for $q > 0$, given by \cite{Wang2005},
\begin{equation*}
\mathrm{vol}(B_p^d)=\frac{\bigl(2\Gamma(1+\frac{1}{p})\bigr)^d}{\Gamma(1+\frac{d}{p})}.
\end{equation*}

Now, taking $f(u) = \e{-\alpha \|u\|^p_p}$ and using the fact that $f(rz) = \e{-\alpha r^p}$ for $z \in \partial B_p^d$ (by definition of the ball), and knowing that $\mu_p$ is a probability measure on $\partial B_p^ d$, results in,
\allowdisplaybreaks
\begin{equation*}
\begin{aligned}
\int_{\mathbb{R}^d} \e{- \alpha \|u\|_p^p} \d{u} =& d\,\mathrm{vol}(B_p^d)\int_0^\infty r^{d-1} \e{-\alpha r^p} \left(\int_{\partial B_p^d} \d{\mu_p(z)}\right)\d{r}\\
    =& d\,\mathrm{vol}(B_p^d) \int_0^\infty e^{- t} \Big(\frac{t}{\alpha }\Big)^{\frac{d-1}{p}} \frac{1}{p \alpha^{1/p}} t^{\frac{1}{p} -1}\d{t} \\ 
    =& d\,\mathrm{vol}(B_p^d) \int_0^\infty e^{- t} \Big(\frac{1}{\alpha }\Big)^{\frac{d}{p}} \frac{1}{p} t^{\frac{d}{p} -1}\d{t} \\
    =& d\,\mathrm{vol}(B_p^d) \frac{\Gamma(d/p)}{p\alpha^{d/p}} 
\end{aligned}
\end{equation*}
where in the second equality we used the change of variable $t = \alpha r^p$ and in the last equality we used the definition of the gamma function, $\Gamma(d) = \int_{0}^{\infty} x^{d-1} e^{-x} \d{x}$. 
\end{proof}

This result together with the following dual formulation of the A-MEP problem, plays a central role in characterizing the Wasserstein-constrained solution.

\begin{proposition}[Dual formulation]
\label{proposition:DRE}
Consider $\mathcal{X} \subseteq \mathbb{R}^d$, then, given $\P{X}{} \in \mathcal{P}_p(\mathcal{X})$, the optimization problem in \cref{equation:average_constrained_entropy} has the following dual formulation,
\begin{customopti*}|s|
    {min}
    {\alpha \geq 0}
    {\alpha \zeta^p + \E{X}{\log \int_\mathcal{X}{\e{- \alpha \|x-v\|^p_p}}\d{v}}}
    {}
    {},
\end{customopti*}
and Radon-Nikodym derivative, 
\begin{equation*}
   \frac{d\P{V|X}{}}{d\P{V}{}} = \frac{\frac{\e{-\alpha \|x-v\|^p_p}}{p(v)}}{\E{v}{\frac{\e{-\alpha \|x-v\|^p_p}}{p(v)}}}.
\end{equation*}
Furthermore, if Slater's condition holds, then strong duality holds and the KKT conditions are necessary and sufficient for optimality. In particular,
\begin{equation*}
\alpha^\star \ge 0, \quad
\Eonly{X}\Eonly{}_{V|X}^\star{\|x-v\|_p^p} \le \zeta^p, \quad
\alpha^\star\left(\Eonly{X}\Eonly{}_{V|X}^\star \left[ \|x-v\|_p^p \right]-\zeta^p\right)=0.
\end{equation*}
\end{proposition}

\begin{proof}
\label{proof:proposition:DRE}
From \cref{lemma:equivalent_formulation} solving \cref{equation:average_constrained_entropy} is equivalent to solving,

\begin{customopti}|s|
    {sup}
    {\P{V|X}{} \in \mathcal{K}_p(\mathcal{X})}
    {\E{X}{\E{V|X}{-\log p(v)} - \KL{\P{V|X}{}}{\P{V}{}}}}
    {}
    {\label{equation:proposition:AMEP_alternative_formulation}}
    \addConstraint{\Eonly{X} \E{V|X}{\|x-v\|_p^p} \leq \zeta^p.}
\end{customopti}
Now, taking the Lagrangian dual,
\begin{customopti*}|s|
    {sup}
    {\begin{aligned}
        \P{V|X}{} \\
        \in \mathcal{K}_p(\mathcal{X})
    \end{aligned}}
    {\Opt{min}
         {\alpha \geq 0}
         {\E{X}{\E{V|X}{-\log p(v)} - \KL{\P{V|X}{}}{\P{V}{}}} - \alpha \left(\E{X}{\E{V|X}{\|x-v\|^p_p}} - \zeta^p \right)}
         {}
    }
    {}
    {},
\end{customopti*}
and grouping the terms,
\begin{customopti}|s|
    {sup}
    {\P{V|X}{} \in \mathcal{K}_p(\mathcal{X})}
    {\underset{\displaystyle \alpha \geq 0}{\min}\quad \E{X}{ \E{V|X}{-\log p(v) - \alpha \|x-v\|^p_p} - \KL{\P{V|X}{}}{\P{V}{}} - \alpha \zeta^p}}
    {}
    {\label{equation:AMEP:proof_max_min}}.
\end{customopti}
At this point, all we can say is that weak duality holds, that is, \cref{equation:AMEP:proof_max_min} $\leq$ \cref{equation:AMEP:proof_min_max}, 
\begin{customopti}|s|
    {sup}
    {\P{V|X}{} \in \mathcal{K}_p(\mathcal{X})}
    {\E{X}{ \E{V|X}{-\log p(v) - \alpha \|x-v\|^p_p} - \KL{\P{V|X}{}}{\P{V}{}} + \alpha \zeta^p}}
    {}
    {\underset{\displaystyle \alpha \geq 0}{\min}}
    {\label{equation:AMEP:proof_min_max}}.
\end{customopti}
In the following steps we will show that \cref{equation:AMEP:proof_min_max} actually has an analytical solution. First, note that, as $\P{X}{}$ is independent of $\P{V|X}{}$, we can look at the expected value of the $\sup$,
\begin{customopti*}|s|
    {sup}
    {\P{V|X}{} \in \mathcal{K}_p(\mathcal{X})}
    {\left. \E{V|X}{-\log p(v) - \alpha \|x-v\|^p_p} - \KL{\P{V|X}{}}{\P{V}{}}  + \alpha \zeta^p \right]}
    {}
    {\underset{\displaystyle \alpha \geq 0}{\min} \Eonly{X} \left[ \right.}.
\end{customopti*}

Now take $\mathcal{X} = \mathbb{R}^d$ and consider,
\begin{equation*}
\begin{aligned}
    \E{V}{\e{-\log p(v) - \alpha \|x-v\|^p_p}} &= \E{V}{\frac{\e{- \alpha \|x-v\|^p_p}}{p(v)}} \\ 
    &= \int_\mathcal{X} \frac{\e{- \alpha \|x-v\|^p_p}}{p(v)} \dP{V}{} \\ 
    &= \int_\mathcal{X} \frac{\e{- \alpha \|x-v\|^p_p}}{p(v)} p(v) \d{v} \\ 
    &= \int_\mathcal{X} \e{- \alpha \|x-v\|^p_p} \d{v}  
\end{aligned}
\end{equation*}
By \cref{proposition:lq_radial_integral}, we know that the above is $< \infty$ if $\alpha >0$. Hence, we can apply Donsker-Varadahn/Gibbs variational principle (\cref{lemma:DV}), and rewrite the inner $\sup$,
\begin{customopti*}|s|
    {min}
    {\alpha > 0}
    {\alpha \zeta^p + \E{X}{\log \int_\mathcal{X}{\e{- \alpha \|x-v\|^p_p}}\d{v}}}
    {}
    {}.
\end{customopti*}

For the Radon-Nikodym derivative, we can use once more Donsker-Varadahn/Gibbs (\cref{lemma:DV}),
\begin{equation*}
   \frac{d \mathbb{P}^\star_{V|X}}{d\P{V}{}} = \frac{\frac{\e{-\alpha \|x-v\|^p_p}}{p(v)}}{\E{V}{\frac{\e{-\alpha \|x-v\|^p_p}}{p(v)}}}.
\end{equation*}

Since the problem is convex, to guarantee strong duality it suffices to verify Slater’s condition, i.e., the existence of a strictly feasible point (an interior point in the relative sense) satisfying the inequality constraints with strict inequality. Under Slater’s condition, the KKT conditions are necessary and sufficient for optimality \cite{Book_Boyd2004}, strong duality holds and optimal Lagrange multipliers exist. Consequently, the KKT conditions are necessary and sufficient for optimality and verifying KKT for the primal-dual optimal, $(\alpha^\star, \P{V|X}{}^\star)$, certifies primal optimality and a zero duality gap \cite{Book_Boyd2004}.
\end{proof}

\subsubsection{Proof of \cref{proposition:AMEP}}
\label{proof:proposition:AMEP}
We finally have all the necessary tools for the proof of \cref{proposition:AMEP},
\begin{proof}
From \cref{proposition:DRE}, we know that \cref{equation:average_constrained_entropy} has a dual representation,
\begin{customopti}|s|
    {min}
    {\alpha > 0}
    {\alpha \zeta^p + \E{X}{\log \int_{\mathbb{R}^d}{\e{- \alpha \|x-v\|^p_p}}\d{v}}}
    {}
    {\label{eq:proposition:AMEPc_afterDV}}.
\end{customopti}

Now, using the result of \cref{proposition:lq_radial_integral} for the integral and \cref{eq:proposition:AMEPc_afterDV} becomes,
\begin{customopti*}
    {min}
    {\alpha > 0}
    {\alpha\zeta^p + \log \Big( \frac{C_{p,d}}{\alpha^{d/p}} \Big)}
    {}
    {},
\end{customopti*}
where $C_{p,d} = \frac{d}{p}\,\mathrm{vol}(B_p^d) \Gamma(d/p)$.

At last, this optimization is solvable and has optimal argument $\alpha^\star = \frac{d}{p\zeta^p}$, since, for $f(\alpha) = \alpha\zeta^p +\log C_{p,d} -\log\alpha^{d/p}$,
\begin{equation*}
\begin{aligned}
    \frac{\partial}{\partial \alpha} f(\alpha) &= \zeta^p - \alpha^{-d/p} \frac{d}{p} \alpha^{d/p-1}\\
    &= \zeta^p - \frac{d}{p} \alpha^{-1}
\end{aligned}
\end{equation*}
and optimal value $\AMEPOptVal$. Leading to the result, by duality,
\begin{customopti*}|s|
    {sup}
    {\P{V|X}{} \in \mathcal{P}}
    {\E{V|X}{-\log p(v)} - \KL{\P{V|X}{}}{\P{V}{}}}
    {}
    {\AMEPOptVal \geq }
    \addConstraint{\Eonly{X} \E{V|X}{\|x-v\|_p^p} \leq \zeta^p.}
\end{customopti*}

For the optimal measure, we know the Radon-Nikodym derivative from \cref{proposition:DRE},
\begin{equation*}
   \frac{d \mathbb{P}^\star_{V|X}}{d\P{V}{}} = \frac{\frac{\e{-\alpha \|x-v\|^p_p}}{p(v)}}{\E{V}{\frac{\e{-\alpha \|x-v\|^p_p}}{p(v)}}}.
\end{equation*}
This tells us that,
\allowdisplaybreaks
\begin{equation*}
\begin{aligned}
   \mathbb{P}^\star_{V|X}(C, \alpha) &= \int_C \frac{\frac{\e{-\alpha \|x-v\|^p_p}}{p(v)}}{\E{V}{\frac{\e{-\alpha \|x-v\|^p_p}}{p(v)}}} \dP{V}{} \\
                  &= \int_C \frac{\frac{\e{-\alpha \|x-v\|^p_p}}{p(v)}}{\E{V}{\frac{\e{-\alpha \|x-v\|^p_p }}{p(v)}}} p(v) \d{v} \\
                  &= \int_C \frac{\e{-\alpha \|x-v\|^p_p }}{\E{V}{\frac{\e{-\alpha \|x-v\|^p_p }}{p(v)}}} \d{v} \\
                  &= \int_C \frac{\e{-\alpha \|x-v\|^p_p }}{\int_{\mathbb{R}^d} \e{-\alpha \|x-v\|^p_p} \d{v}} \d{v}.
\end{aligned}
\end{equation*}

Resulting in an optimal argument,
\begin{equation*}
\mathbb{P}^\star_{V|X}(C) = \AMEPOptArg.
\end{equation*}

Now, to show strong duality we need to show Slater and KKT conditions as per \cref{proposition:DRE}:

\paragraph{1. Slater condition} We need to check that there exists a strictly feasible point for \cref{equation:proposition:AMEP_alternative_formulation}. Assuming $\zeta >0$, take $U$ to be a uniform distribution on $[-\delta,\delta]^d$ and independent of $X$. Define $V=X+U$, then $\E{X,V}{\|x-v|^p_p} = E[\|u\|^p_p] \leq C(q) \delta^p$ (with $C(q)$ being a scaling factor to account for the norm used). Choosing $\delta$ small enough such that $C(q) \delta^p < \zeta^p$ guarantees that Slater's condition is satisfied.

\paragraph{2. KKT conditions} 
\begin{itemize}
  \item Dual feasibility: $\alpha^\star\ge 0$. 
  \item Primal feasibility: $\int_{\mathcal{X}} \int_{\mathcal{X}} \|x-v\|^p_p \dP{V|X}{}^\star(v) \dP{X}{x} \leq \zeta^p$.
  \item Complementary slackness: $\alpha^\star\big(\zeta^p- \int_{\mathcal{X}} \int_{\mathcal{X}} \|x-v\|^p_p \dP{V|X}{}^\star(v) \dP{X}{x} \big)=0$.
  \item Stationarity: $\P{V|X}{}^\star$ maximizes the Lagrangian at $\alpha^\star$. 
\end{itemize}

Dual feasibility holds immediately by inspecting $\alpha^\star = \frac{d}{p \zeta^p}$. Similarly, stationarity holds by construction, since $\P{V|X}{}^\star$ maximizes the Lagrangian for $\alpha^\star$ via Donsker-Varadahn/Gibbs result. Since $\alpha^\star >0$, it suffices to show $\int_{\mathcal{X}} \int_{\mathcal{X}} \|x-v\|^p_p \dP{V|X}{}^\star(v) \dP{X}{x} = \zeta^p$ for the primal feasibility and complementary slackness,
\begin{equation*}
    \begin{aligned}
        \allowdisplaybreaks
        &= \int_{\mathcal{X}} \int_{\mathcal{X}} \frac{\|x-v\|^p_p \e{-\alpha^\star \|x-v\|^p_p}}{ \int_{\mathcal{X}} \e{-\alpha^\star \|x-v\|^p_p} \d{v}} \d{v} \dP{X}{x} \\
        &= \int_{\mathcal{X}} \frac{d\mathrm{vol}(B_p^d) \int_0^\infty  r^p r^{d-1} \e{-\alpha^\star r^p} \d{r}}{d\mathrm{vol}(B_p^d) \int_0^\infty r^{d-1}\e{-\alpha^\star r^p}} \d{r} \dP{X}{x} \\
        &= \int_{\mathcal{X}} \frac{\int_0^\infty r^p r^{d-1} \e{-\alpha^\star r^p} \d{r} }{\int_0^\infty r^{d-1}\e{-\alpha^\star r^p} \d{r}} \dP{X}{x} \\
        &= \int_{\mathcal{X}} \frac{\int_0^\infty r^{p+d-1} \e{-\alpha^\star r^p} \d{r} }{\int_0^\infty r^{d-1}\e{-\alpha^\star r^p} \d{r}} \dP{X}{x} \\
        &= \int_{\mathcal{X}} \frac{\frac{\Gamma(1+d/p)}{p {\alpha^\star}^{1+d/p}}}{\frac{\Gamma(d/p)}{p {\alpha^\star}^{d/p}}} \dP{X}{x} \\
        &= \int_{\mathcal{X}} \frac{\Gamma(1+d/p)}{\Gamma(d/p)\alpha^\star} \dP{X}{x} \\
        &= \int_{\mathcal{X}} \frac{d/p\Gamma(d/p)}{\Gamma(d/p)\alpha^\star} \dP{X}{x} \\
        &= \frac{d}{p\alpha^\star}\\
        &= \zeta^p 
    \end{aligned}
\end{equation*}
where we used the property of the Gamma function, $\Gamma(1+z)= z\Gamma(z)$.

Consequently, strong duality holds, leading to,
\begin{customopti*}|s|
    {sup}
    {\P{V|X}{} \in \mathcal{P}}
    {\DiffEnt{\P{X,V}{}}{V|X}}
    {}
    {\AMEPOptVal =}
    \addConstraint{\Eonly{X} \E{V|X}{\|x-v\|_p^p} \leq \zeta^p},
\end{customopti*}
which concludes the proof.
\end{proof}

\subsubsection{Proof of \cref{corollary:AMEP2}}
\label{proof:corollary:AMEP2}
\begin{proof}
    Setting $p=2$ in \cref{proposition:AMEP} gives us,
    \begin{equation*}
\frac{d}{2} + \log \Big(\sqrt{\frac{2\pi}{d}} \zeta\Big)^{d},
\end{equation*}
for the optimal solution and
\begin{equation*}
\mathbb{P}^\star_{V|X}(C) =  \frac{1}{\Big(\sqrt{\frac{2 \pi}{d}}\zeta \Big)^d} \int_C \e{-\frac{d}{2\zeta^2} \|x-v\|_2^2} \d{v},
\end{equation*}
for the optimal argument.
Since $x$ is given and both $x,v\in \mathbb{R}^d$, we can write the above as,
\begin{equation*}
\frac{1}{\sqrt{(2 \pi)^d \mathrm{det}(\frac{\zeta^2}{d}\Id{d})}} \int_C \e{-\frac{1}{2} \transpose{(v-x)} (\frac{\zeta^2}{d}\Id{d})^{-1} (v-x)} \d{v},
\end{equation*}
which is a  Gaussian kernel, $\mathcal N\!\left(x,\frac{\zeta^2}{d}\Id{d}\right)$.
\end{proof}

\subsubsection{Proof of \cref{corollary:AMEP1}}
\label{proof:corollary:AMEP1}
\begin{proof}
    Setting $p=1$ in \cref{proposition:AMEP} gives us,
\begin{equation*}
d + \log \Big( \frac{2 \zeta}{d}\Big)^d,
\end{equation*}
for the optimal solution and 
\begin{equation*}
\mathbb{P}^\star_{V|X}(C) = \Big(\frac{d}{2 \zeta}\Big)^d \int_C \e{-\frac{d}{\zeta} \|x-v\|_1} \d{v}.
\end{equation*}
for the optimal argument, which is equivalent to writing, 
\[
\begin{aligned}
\mathbb{P}^\star_{V|X}(C) &= \Big(\frac{d}{2 \zeta}\Big)^d \int_C \e{-\frac{d}{\zeta} \sum_{i=1}^d \|x_i-v_i\|_1} \d{v}_1 \ldots \d{v}_d\\
 &= \Big(\frac{d}{2 \zeta}\Big)^d \int_C \prod_{i=1}^d \e{-\frac{d}{\zeta} \|x_i-v_i\|_1} \d{v}_1 \ldots \d{v}_d,\\
 &=  \int_C \prod_{i=1}^d \frac{d}{2 \zeta} \e{-\frac{d}{\zeta} \|x_i-v_i\|_1} \d{v}_1 \ldots \d{v}_d.
\end{aligned}
\]
Since the above holds for every Borel set $C\subseteq\mathbb{R}^d$, the optimal
conditional law $P^\star_{V|X}$ is absolutely continuous w.r.t. Lebesgue measure with density
\[
p^\star(v|x)=\prod_{i=1}^d \frac{d}{2\zeta}\exp\!\Big(-\frac{d}{\zeta}\|x_i-v_i\|_1\Big).
\]
This is the joint density of independent coordinates $v_i|x \stackrel{\mathrm{i.i.d.}}{\sim} \mathrm{Laplace}(x_i,\zeta/d)$, i.e.\ $V| x \sim \mathrm{Laplace}_d(x,\zeta/d)$.
\end{proof}

\subsubsection{Proof of \cref{corollary:AMEPinfty}}
\label{proof:corollary:AMEPinfty}
\begin{proof}
The constraint restricts the support of $P_{V|X=x}$ to
\[
D_\zeta(x):=\bigl(x+[-\zeta,\zeta]^d\bigr),
\]
Among all probability measures supported on a measurable set, the uniform distribution uniquely maximizes differential entropy. Hence
\[
\P{V|X=x}{}^\star = \mathrm{Unif}(D_\zeta(x)),
\]
with optimal value
\[
\DiffEnt{\P{X,V}{}^\star}{V|X=x}=\log|D_\zeta(x)|.
\]

If $x\in[\zeta,1-\zeta]^d$, then
\[
D_\zeta(x)=x+[-\zeta,\zeta]^d,
\]
whose Lebesgue measure is $(2\zeta)^d$. Therefore
\[
\DiffEnt{\P{X,V}{}^\star}{V|X=x} = \log(2\zeta)^d = d\log(2\zeta).
\]

\end{proof}

\subsection{Proof of W-MEP results}

Prior to the proceeding with the proof we first show a simple statement which will be used in this and subsequent proofs,

\begin{lemma}
\label{lemma:sup_measure_to_pointwise_swap}
Let $(\mathcal{X},\mathcal{A})$ be a measurable space and let $f:\mathcal{X}\to\mathbb{R}$ be measurable. Then
\[
\Opt{sup}{\P{}{} \in \mathcal{P}(\mathcal{X})}{\int_{\mathcal{X}} f(v) \dP{}{v}}{} =\Opt{sup}{v \in \mathcal{X}}{f(v)}{}
\]
provided that the integrals are well defined.
\end{lemma}

\begin{proof}
Let 
\[
M := \Opt{sup}{v \in \mathcal{X}}{f(v)}{},
\]
then, for every $\P{}{} \in\mathcal{P}(\mathcal{X})$,
\[
\int_{\mathcal{X}} f(v) \dP{}{v} \leq \int_{\mathcal{X}} M \dP{}{v}.
\]
Now, since $\P{}{}$ is a probability measure, then,
\[
\int_{\mathcal{X}} f(v) \dP{}{v} \leq M.
\]
Conversely, for every $v\in\mathcal{X}$, the Dirac measure $\delta_v$ belongs to $\mathcal{P}(\mathcal{X})$, and
\[
\int_{\mathcal{X}} f(u)\,\delta_v(du) = f(v).
\]
Therefore,
\[
\Opt{sup}{\P{}{}\in\mathcal{P}(\mathcal{X})}{\int_{\mathcal{X}} f(u) \dP{}{u}}{}  \geq f(v) \qquad \text{for every } v\in\mathcal{X}.
\]
Taking the supremum over $v\in\mathcal{X}$ yields
\[
\Opt{sup}{\P{}{}\in\mathcal{P}(\mathcal{X})}{\int_{\mathcal{X}} f(u) \dP{}{u}}{}  \geq \Opt{sup}{v \in \mathcal{X}}{f(v)}{}.
\]
Combining the two inequalities proves the result.
\end{proof}

\subsubsection{Proof of \cref{proposition:WMEP}}
\label{proof:proposition:WMEP}
\begin{proof}
    
From \cref{lemma:equivalent_formulation} solving \cref{equation:wasserstein_constrained_entropy} is equivalent to solving,
\begin{customopti*}|s|
    {sup}
    {\P{V|X}{} \in \mathcal{K}_p(\mathcal{X})}
    {\E{X}{\E{\P{V|X}{}}{-\log p(v)} - \KL{\P{V|X}{}}{\P{V}{}}}}
    {}
    {}
    \addConstraint{W_p(\P{X}{}, \P{V}{}) \leq \zeta}.
\end{customopti*}
Note that it is equivalent to (by the definition of the Wasserstein metric),

\begin{customopti*}|s|
    {sup}
    {\P{V|X}{} \in \mathcal{K}_p(\mathcal{X})}
    {\E{X}{\E{\P{V|X}{}}{-\log p(v)} - \KL{\P{V|X}{}}{\P{V}{}}}}
    {}
    {}
    \addConstraint{\Opt{inf}{\PQ{}{} \in \Pi(\P{V}{}, \P{X}{})}{\int_{\mathcal{X} \times \mathcal{X} } \|x-v\|^p_p \dPQ{}{x,v}}{} \leq \zeta^p },
\end{customopti*}
written differently, 
\begin{customopti*}|s|
    {sup}
    {\begin{aligned}
    \P{V|X}{} &\in \mathcal{K}_p(\mathcal{X}) \\
    \pi \in & \Pi(\P{X}{}, \P{V}{})
    \end{aligned}}
    {\DiffEnt{\P{X,V}{}}{V|X}}
    {}
    {}
    \addConstraint{\int_{\mathcal{X} \times \mathcal{X} } \|x-v\|^p_p \dPQ{}{x,v}  \leq \zeta^p }.
\end{customopti*}
Since $(\mathcal X,\mathcal B(\mathcal X))$ is a Borel space, every coupling
$\PQ{}{}\in\Pi(\P{X}{},\P{V}{})$ admits a disintegration
\[
\dPQ{}{x,v} = \dP{X}{x}\,\dQ{V|X}{v},
\]
where $\Q{V|X}{}$ is the transport kernel.
In contrast to the average transport formulation, this kernel need not coincide with the conditional law $\P{V|X}{}$ appearing in the entropy objective. Indeed, the Wasserstein constraint depends only on the induced marginal distribution $\P{V}{}$, whereas the entropy depends on $\P{V|X}{}$. Consequently, the optimization naturally involves two kernels, $\P{V|X}{}$ and $\Q{V|X}{}$, which are linked only through the consistency condition
\[
\int_{\mathcal X}\P{V|X}{C}\dP{X}{x} = \int_{\mathcal X}\Q{V|X}{C}\dP{X}{x}, \qquad \forall C \in \mathcal{B}(\mathcal{X}).
\]
Therefore, optimizing over couplings $\pi \in \Pi(\P{X}{}, \P{V}{})$ is equivalent to optimizing over transport kernels $\Q{V|X}{}$ satisfying the marginal consistency condition. Hence, the problem can be reformulated over the pair $\P{V|X}{},\Q{V|X}{}$ as, 
\begin{customopti*}|s|
    {sup}
    {\P{V|X}{}, \Q{V|X}{} \in \mathcal{K}_p(\mathcal{X})}
    {\E{X}{\E{\P{V|X}{}}{-\log p(v)} - \KL{\P{V|X}{}}{\P{V}{}}}}
    {}
    {}
    \addConstraint{\int_{\mathcal{X}} \int_{\mathcal{X}} \|x-v\|^p_p \dQ{V|X}{v} \dP{X}{x}  \leq \zeta^p }
    \addConstraint{\int_{\mathcal X}\P{V|X}{C}\dP{X}{x} = \int_{\mathcal X}\Q{V|X}{C}\dP{X}{x}, \quad \forall C \in \mathcal{B}(\mathcal{X})}.
\end{customopti*}
which is equivalent to the following (expressing the integrals as expectations),
\begin{customopti}|s|
    {sup}
    {\P{V|X}{}, \Q{V|X}{} \in \mathcal{K}_p(\mathcal{X})}
    {\E{X}{\E{\P{V|X}{}}{-\log p(v)} - \KL{\P{V|X}{}}{\P{V}{}}}}
    {}
    {}
    \addConstraint{\E{X}{\E{\Q{V|X}{}}{\|x-v\|^p_p}} \leq \zeta^p }
    \addConstraint{\int_{\mathcal X}\P{V|X}{C}\dP{X}{x} = \int_{\mathcal X}\Q{V|X}{C}\dP{X}{x}, }
    \addConstraint{\qquad \qquad \forall C \in \mathcal{B}(\mathcal{X})}.
\end{customopti}
The consistency constraint is an equality between probability measures rather than a scalar constraint. Following the standard convex-analytic duality for optimization over measures (see, e.g., \cite{Villani2008,Santambrogio2015}), its Lagrange multiplier is therefore a measurable potential $\varphi:\mathcal{X} \rightarrow \mathbb{R}$, which acts on measures through the inner product 
\[
\langle\varphi,\mu\rangle = \int_{\mathcal{X}}\varphi(v)\d{\mu(v)}.
\]
Introducing the multiplier $\alpha\ge0$ for the transport budget and the
potential $\varphi$ for the marginal consistency constraint yields the
dual Lagrangian problem,
\begin{customopti*}|s|
    {sup}
    {\P{V|X}{}, \Q{V|X}{} \in \mathcal{K}_p(\mathcal{X})}
    {\Opt{min}
         {\alpha \geq 0, \varphi}
         {\left\{ \begin{aligned}
            &\E{X}{\E{\P{V|X}{}}{-\log p(v)} - \KL{\P{V|X}{}}{\P{V}{}}} \\
            &\quad - \alpha \E{X}{\E{\Q{V|X}{}}{\|x-v\|^p_p}} + \alpha \zeta^p \\
            &\quad - \E{X}{\E{\P{V|X}{}}{\varphi(v)}} + \E{X}{\E{\Q{V|X}{}}{ \varphi(v)}}
         \end{aligned} \right\}}
         {}
    }
    {}
    {},
\end{customopti*}
and grouping the terms,
\begin{customopti}|s|
    {sup}
    {\begin{aligned}
    \P{V|X}{}, \Q{V|X}{} \\
    \in \mathcal{K}_p(\mathcal{X})
    \end{aligned}}
    {\Opt{inf}
         {\alpha \geq 0, \varphi}
         {\left\{ \begin{aligned}
         &\E{X}{\E{\P{V|X}{}}{-\varphi(v) -\log p(v) } - \KL{\P{V|X}{}}{\P{V}{}}} \\
         &+ \E{X}{\E{\Q{V|X}{}}{ \varphi(v) - \alpha \|x-v\|^p_p}} + \alpha \zeta^p
         \end{aligned}\right\}
         }
         {}
         {}
    }
    {\label{equation:WMEP:proof_max_min}}
    {}.
\end{customopti}
At this point, all we can say is that weak duality holds, that is, \cref{equation:WMEP:proof_max_min} is less or equal to \cref{equation:WMEP:proof_min_max}, 
\begin{customopti}|s|
    {inf}
    {\alpha \geq 0, \varphi}
    {\Opt{sup}
        {\begin{aligned}
            \P{V|X}{}, \Q{V|X}{} \\
            \in \mathcal{K}_p(\mathcal{X})
        \end{aligned}}
         {\left\{ \begin{aligned}
         &\E{X}{\E{\P{V|X}{}}{-\varphi(v) -\log p(v) } - \KL{\P{V|X}{}}{\P{V}{}}} \\
         &+ \E{X}{\E{\Q{V|X}{}}{ \varphi(v) - \alpha \|x-v\|^p_p}} + \alpha \zeta^p
         \end{aligned}\right\}
         }
         {}
         {}
    }
    {\label{equation:WMEP:proof_min_max}}
    {}.
\end{customopti}
The objective is separable in $\P{V|X}{}$ and $\Q{V|X}{}$, that is, we can optimize over each argument separately.
\paragraph{Entropy Kernel, $\P{V|X}{}$}
\begin{customopti*}|s|
    {sup}
    {\P{V|X}{} \in \mathcal{K}_p(\mathcal{X})}
    {\E{X}{\E{\P{V|X}{}}{-\varphi(v) -\log p(v) } - \KL{\P{V|X}{}}{\P{V}{}}}} 
    {}
    {}.
\end{customopti*}
Note that, as $\P{X}{}$ is independent of $\P{V|X}{}$, we can look at the expected value of the $\sup$,
\begin{customopti*}|s|
    {sup}
    {\P{V|X}{} \in \mathcal{K}_p(\mathcal{X})}
    {\E{\P{V|X}{}}{-\varphi(v) -\log p(v) } - \KL{\P{V|X}{}}{\P{V}{}}  \Big]}
    {}
    {\Eonly{X} \Big[}.
\end{customopti*}

Now take $\mathcal{X} = \mathbb{R}^d$ and consider,
\begin{equation*}
\begin{aligned}
    \E{V}{\e{-\varphi(v) -\log p(v)}} &= \int_\mathcal{X} \frac{\e{-\varphi(v)}}{p(v)} \dP{V}{} \\ 
    &= \int_\mathcal{X} \frac{\e{-\varphi(v)}}{p(v)} p(v) \d{v} \\ 
    &= \int_\mathcal{X} \e{-\varphi(v)} \d{v}  
\end{aligned}
\end{equation*}
By the Donsker--Varadhan/Gibbs variational principle, \cref{lemma:DV}, the supremum over $\P{V|X}{}$ is attained at a kernel satisfying, for $\P{X}{}$-a.e. $x$,
\[
\frac{\dP{V|X}{}^{\star}(\cdot|x)}{\dP{V}{}}(v) = \frac{\e{-\varphi(v)-\log p(v)}}{\E{\P{V}{}}{\e{-\varphi(V)-\log p(V)}}}
\]
Since $\dP{V}{}(v)=p(v)\d{v}$, the denominator becomes
\[
\E{V}{\e{-\varphi(V)-\log p(V)}} = \int_{\mathcal X} \e{-\varphi(v)}\d{v}.
\]
Therefore,
\[
\dP{V|X}{}^{\star}(v) = \frac{\e{-\varphi(v)}}{\int_{\mathcal X}\e{-\varphi(u)}\d{u}}\d{v}, \qquad \P{X}{}\text{-a.s.}
\]
In particular, the entropy-maximizing kernel is independent of $x$. Hence the
induced marginal satisfies
\[
p_V^\star(v) = \frac{\e{-\varphi^\star(v)}} {\int_{\mathcal X}\e{-\varphi^\star(u)}\d{u}}.
\]
Consequently,
\[
\P{V|X}{}^{\star}(\cdot)=\P{V}{}^{\star}(\cdot), \qquad \P{X}{}\text{-a.s.},
\]
and the conditional entropy equals the marginal entropy at the optimum.

\paragraph{Transport Kernel, $\Q{V|X}{}$}
\begin{customopti*}|s|
    {sup}
    {\Q{V|X}{} \in \mathcal{K}_p(\mathcal{X})}
    {\E{X}{\E{\Q{V|X}{}}{ \varphi(v) - \alpha \|x-v\|^p_p}}} 
    {}
    {}.
\end{customopti*}
Similarly, we can look at,
\begin{customopti*}|s|
    {sup}
    {\Q{V|X}{} \in \mathcal{K}_p(\mathcal{X})}
    {\E{\Q{V|X}{}}{ \varphi(v) - \alpha \|x-v\|^p_p} \Big]} 
    {}
    {\Eonly{X} \Big[},
\end{customopti*}
and from \cref{lemma:sup_measure_to_pointwise_swap}, this is equivalent to,
\begin{customopti*}|s|
    {sup}
    {v \in \mathcal{X}}
    {\varphi(v) - \alpha \|x-v\|^p_p \Big] } 
    {}
    {\Eonly{X} \Big[ },
\end{customopti*}
which is,
\[
\E{X}{-\varphi^{c_\alpha}(x)},
\]
since, the $c$-transform \cite{Villani2008, Santambrogio2015} for $c_\alpha(x,v) = \alpha \| x - v\|_p^p$ is defined as, 
\[
\varphi^{c_\alpha}(x) = \Opt{inf}{v \in \mathcal{X}}{\alpha \|x-v \|_p^p-\varphi(v)}{}.
\]

Combining these results into \cref{equation:WMEP:proof_min_max}, results in, 
\begin{customopti*}|s|
    {inf}
    {\alpha \geq 0, \varphi}
    {\log \int_\mathcal{X} \e{-\varphi(v)} \d{v} - \E{X}{\varphi^{c_\alpha}(x)} + \alpha \zeta^p}
    {}
    {}.
\end{customopti*}
Since the problem is convex, to guarantee strong duality it suffices to verify Slater’s condition, i.e., the existence of a strictly feasible point (an interior point in the relative sense) satisfying the inequality constraints with strict inequality. Under Slater’s condition, the KKT conditions are necessary and sufficient for optimality \cite{Book_Boyd2004}, strong duality holds and optimal Lagrange multipliers exist. Consequently, the KKT conditions are necessary and sufficient for optimality and verifying KKT for the primal-dual optimal, $(\alpha^\star, \varphi^\star)$, certifies primal optimality and a zero duality gap \cite{Book_Boyd2004}.

\end{proof}

\subsubsection{Proof of \cref{corollary:WMEP1}}
\label{proof:corollary:WMEP1}
\begin{proof}
The result follows directly from \cref{proposition:WMEP}. For $p=1$, let $\varphi$ be $\alpha$-Lipschitz. Then
\[
\varphi(v)-\varphi(x)\leq \alpha\|v-x\|_1,
\]
and hence
\[
\alpha\|x-v\|_1-\varphi(v)\geq -\varphi(x).
\]
Since equality is attained at $v=x$, the $c_\alpha$-transform satisfies
\[
\varphi^{c_\alpha}(x) = \Opt{inf}{v}{\alpha\|x-v\|_1-\varphi(v)}{} = -\varphi(x).
\]
Substituting into the dual formulation of \cref{proposition:WMEP} gives
\[
\Opt{inf}{\alpha\geq0,\\\varphi\in\mathrm{Lip}_\alpha}{\log\int_{\mathbb R^d}\e{-\varphi(v)} \d{v} +\E{X}{\varphi(X)} +\alpha\zeta}{}
\]
For a fixed Lipschitz potential $\varphi$, the smallest admissible $\alpha$ is $\operatorname{Lip}(\varphi)$, which yields the stated formulation. The expression for the optimal marginal and the independence of the optimal conditional kernel from $X$ follow directly from \cref{proposition:WMEP}.
\end{proof}

\end{document}